\documentclass[VANCOUVER,Times1COL]{WileyNJDv5} %STIX1COL,STIX2COL,STIXSMALL
\articletype{Article Type}%

\received{Date Month Year}
\revised{Date Month Year}
\accepted{Date Month Year}
\journal{Journal}
\volume{00}
\copyyear{2025}
\startpage{1}

\usepackage{multirow} % In ACM package list
\usepackage[font=normal]{subcaption} % In ACM package list
\usepackage{geometry} % In ACM package list
\usepackage{setspace} % In ACM package list

\usepackage{accents} % Not in ACM package list
\usepackage{trimclip} % Not in ACM package list

\usepackage{arydshln}

\graphicspath{{images/}}

\begin{document}

%\renewcommand{\todo}[1]{} % Note: This can be used to hide todo comments

%=====================================
%      my definitions 
%=====================================
\newcommand{\newtilde}{\scalebox{.5}[.4]{\trimbox{0pt 1.8pt}{$\sim$}}}
\newcommand{\newtildeb}{\scalebox{.5}[.4]{\trimbox{0pt 1.4pt}{$\sim$}}}
\newcommand{\M}[1]{\boldsymbol{#1}}

\newcommand{\Mt}[1]{{ \bf #1}}

\newcommand{\bld}[1]{\boldsymbol{#1}}
\newcommand{\inv}[1]{#1^{\scriptscriptstyle -1}}
%\stackMath

% I had a little trouble picking the tilde symbol for the approximated matrices 
% in a way that was both pretty, portable, and required minimal packages. 
% In particular the double tilde is tricky.

%\newcommand{\til}[1]{\tilde{#1}{}}
%\newcommand{\til}[1]{\overset{ \scriptscriptstyle \sim}{#1}}
%\newcommand{\til}[1]{\accentset{ \scriptscriptstyle \sim}{#1}}
%\newcommand{\til}[1]{\stretchedtilde{#1}}
%\newcommand{\til}[1]{\stackengine{.5pt}{#1}{\scriptscriptstyle \sim}{O}{c}{F}{F}{S}}
%\newcommand{\til}[1]{\trlap[.5pt]{\scriptscriptstyle \sim}{#1}}
%\newcommand{\dtil}[1]{\tilde{\tilde{#1}}{}}
%\newcommand{\til}[1]{\stackon[.5pt]{#1}{\scriptscriptstyle \hspace{0pt} \sim}{}}
%\newcommand{\dtil}[1]{\stackon[-.5pt]{\stackon[.5pt]{#1}{\scriptscriptstyle \hspace{0pt} \sim}}{\scriptscriptstyle \hspace{0pt} \sim}{}}

\newcommand{\til}[1]{\accentset{ {\newtilde}}{#1}}
\newcommand{\dtil}[1]{\accentset{ {\newtildeb}}{\til{#1}}}

\newcommand{\spaceh}{\hspace{.5cm}}
\newcommand{\vertcells}[3]{\multirow{#1}{*}{\rotatebox[origin=c]{#2}{\parbox[c]{.33cm*#1}{#3}}}}
\newcommand{\hcells}[3]{\multirow{#1}{*}{\rotatebox[origin=c]{#2}{\parbox[c]{1cm}{#3}}}}
\newcommand{\emin}[1]{10^{-#1}}
\newcommand{\etothe}[1]{10^{#1}}
\newcommand{\Vmat}{\hat{T}}
\newcommand{\vmat}{T}
\newcommand{\Wmat}{\hat{W}}
\newcommand{\wmat}{W}
\newcommand{\Bmat}{\hat{B}}
\newcommand{\bmat}{B}
\newcommand{\Cmat}{\hat{C}}
\newcommand{\cmat}{C}
\newcommand{\red}{r}
\newcommand{\bw}{k}
\newlength{\commentWidth}
\newcommand{\atcp}[1]{\tcp*[r]{\makebox[\commentWidth]{#1\hfill}}}
\newcommand{\optpagebreak}{}
\newcommand{\typicalsize}{\scriptstyle}
\newcommand{\smallmatsize}{\scriptscriptstyle}
\newcommand{\typicaldims}{1pt,5cm,5cm}
\newcommand{\eqn}{Eq.}
\newcommand{\Eqn}{Eq.}
\newcommand{\eqns}{Eqs.}
\newcommand{\Eqns}{Eqs.}
\newcommand{\Fig}{Fig.}
\newcommand{\fig}{Fig.}
\newcommand{\Figs}{Figs.}
\newcommand{\figs}{Figs.}
\newcommand{\spike}{SPIKE}
\newcommand{\Spike}{SPIKE}
\newcommand{\smul}{\hspace{.01cm}}
\newcommand{\prettycaption}{\captionsetup[subfigure]{width=0.8\linewidth,justification=raggedright}}
\newcommand{\algstub}{LR}
\newcommand{\algname}{\algstub-\spike{}}
\newcommand{\Algname}{\algstub-\Spike{}}
\newcommand{\mco}[1]{\multicolumn{2}{c|}{#1}}
\newcommand{\mcol}[1]{\multicolumn{2}{c||}{#1}}
\newcommand{\mct}[2]{$#1$  &\hspace{-1em} $\left(#2\right)$ }
\newcommand{\mctl}[2]{$#1$ &\hspace{-1em} $\left(#2\right)$ }
\newcommand{\tworow}[1]{\multirow{2}{*}{#1}}
\newcommand{\mcttworow}[2]{\tworow{$#1$}  &\hspace{-1em} \tworow{$\left(#2\right)$} }

\makeatletter
\newcommand{\newdecl}[2]{\csgdef{decl@#1}{#2}}% Creates a declaration
\newcommand{\csvdel}{}% Delimiter used in CSV representation
\newcommand{\newusedecl}[2][,]{% Use a declaration
  \renewcommand{\csvdel}{\renewcommand{\csvdel}{#1\,}}% Delay \csvdel one cycle.
  \csname decl@#2\endcsname(\checknextarg}
\newcommand{\checknextarg}{\@ifnextchar\bgroup{\gobblenext}{}}% Check if another "argument" exists
\newcommand{\gobblenext}[1]{\csvdel#1\@ifnextchar\bgroup{\gobblenext}{)}}% Gobble next "argument"
\makeatother

\newcommand{\tempA}{g}
\newcommand{\tempB}{h}

\newcommand{\fillerspace}{}
%\newcommand{\fillerspace}
%{
%\begin{itemize}
%\item
%\item \hrulefill
%\item
%\item \hrulefill
%\item
%\item \hrulefill
%\item
%\item \hrulefill
%\item
%\item \hrulefill
%\item
%\item \hrulefill
%\item
%\item \hrulefill
%\item
%\item \hrulefill
%\item
%\item \hrulefill
%\end{itemize}
%}

\newcommand{\mybmat}[3]{%%
\left[ %%
\typicalsize %%
\begin{array}{#1} %%
#3 %%
\end{array} %%
\right] %%
}

\newcommand{\mybmatA}[3]{%%
\renewcommand\arraystretch{1.3}
\left[ %%
\typicalsize %%
%\begin{BMAT}(@,.55cm,.55cm)[0pt]{#1}{#2} %%
\begin{array}{#1} %%
#3 %%
\end{array} %%
%\end{BMAT} %%
\right] %%
}

\newcommand{\myaltbmatA}[3]{%%
\renewcommand\arraystretch{1.3}
\left[ %%
\smallmatsize %%
\begin{array}{#1} %%
#3 %%
\end{array} %%
\right] %%
}

\newcommand{\myaltbmatB}[3]{ %%
\renewcommand\arraystretch{1.3}
\left[ %%
\typicalsize %%
\begin{array}{#1} %%
#3 %%
\end{array} %%
\right] %%
}

\newcommand{\myaltbmatC}[3]{ %%
\left[ %%
\smallmatsize %%
\begin{array}{#1} %%
#3 %%
\end{array} %%
\right] %%
}

\newcommand{\myaltbmatCII}[3]{ %%
\renewcommand\arraystretch{2}
\left( %%
\smallmatsize %%
\begin{array}{#1} %%
#3 %%
\end{array} %%
\right) %%
}

\newcommand{\notmybmat}[3]{%%
\left[ %%
\typicalsize %%
\begin{BMAT}(b,.55cm,.55cm)[0pt]{#1}{#2} %%
#3 %%
\end{BMAT} %%
\right] %%
}

\newcommand{\notmybmatA}[3]{%%
\left[ %%
\typicalsize %%
\begin{BMAT}(@,.55cm,.55cm)[0pt]{#1}{#2} %%
#3 %%
\end{BMAT} %%
\right] %%
}

\newcommand{\notmyaltbmatA}[3]{%%
\left[ %%
\smallmatsize %%
\begin{BMAT}(e,.7cm,.7cm)[0pt]{#1}{#2} %%
#3 %%
\end{BMAT} %%
\right] %%
}

\newcommand{\notmyaltbmatB}[3]{ %%
\left[ %%
\typicalsize %%
\begin{BMAT}(b,.5cm,.4cm)[0pt]{#1}{#2} %%
#3 %%
\end{BMAT} %%
\right] %%
}

\newcommand{\notmyaltbmatC}[3]{ %%
\left[ %%
\smallmatsize %%
\begin{BMAT}(@,.5cm,.5cm)[2pt,10mm,10mm]{#1}{#2} %%
#3 %%
\end{BMAT} %%
\right] %%
}

%\newcommand{\optpagebreak}{}

%=====================================
%     end my definitions 
%=====================================

%\title{Reduced Rank Approximation Framework for the SPIKE Linear System Solver with Applications}

%\title{Reduced Rank SPIKE Framework for solving sparse Linear System}

\title{Low-Rank SPIKE Framework for Solving Large Sparse Linear Systems with Applications}

%\title{Reduced Rank extensions for the SPIKE family of Linear System Solver}

\author[1]{Braegan S. Spring}

\author[1,2]{Eric Polizzi}

\author[3]{Ahmed H. Sameh}

\authormark{B. S. Spring, E. Polizzi, A. H. Sameh}
\titlemark{Low-Rank Framework for SPIKE}

\address[1]{\orgdiv{Electrical and Computer Engineering}, \orgname{University of Massachusetts, Amherst}, \orgaddress{\state{MA}, \country{USA}}}

\address[2]{\orgdiv{Mathematics and Statistics}, \orgname{University of Massachusetts}, \orgaddress{\state{MA}, \country{USA}}}

\address[3]{\orgdiv{Department of Computer Sciences}, \orgname{Purdue University, West-Lafayette}, \orgaddress{\state{IN}, \country{USA}}}

% TODO figure out name
\corres{Corresponding author Braegan Spring. \email{bspring@umass.edu}}

%\presentaddress{This is sample for present address text this is sample for present address text.}

%\fundingInfo{Text}
%\JELinfo{ejlje}

\abstract[Abstract]{
%The \Spike{} family of linear system solvers provides parallelism targeted to banded systems using a block tridiagonal partitioning. 
%For banded matrices this results in a block tridiagional representation, with off-diagonal blocks representing communication. 

%For banded matrices this results in a block tridiagional representation, with off-diagonal blocks representing communication. 

%The \Spike{} family of linear system solvers provides parallelism using a block partitioning along the matrix diagonal. 
%For banded matrices this results in a block tridiagional representation, with off-diagonal blocks representing communication. 

The \spike{} family of linear system solvers provides parallelism using a block tridiagonal partitioning. 
Typically \spike{}-based solvers are applied to banded systems, resulting in structured off-diagonal blocks with non-zeros elements restricted to relatively small submatrices comprising the band of the original matrix. 
In this work, a low-rank SVD based approximation of the off-diagonal blocks is investigated.
This produces a representation which more effectively handles matrices with large, sparse bands.
A set of flexible distributed solvers, the \algname{} variants, are implemented.
There are applicable to a wide range of applications\textemdash from use as a ``black-box'' preconditioner which straightforwardly improves upon the classic Block Jacobi preconditioner, to use as a specialized ``approximate direct solver.'' 
An investigation of the effectiveness of the new preconditioners for a selection of SuiteSparse matrices is performed, particularly focusing on matrices derived from 3D finite element simulations.
In addition, the \spike{} approximate linear system solvers are also paired with the FEAST eigenvalue solver, where they are shown to be particularly effective due to the former's rapid convergence, and the latter's acceptance of loose linear system solver convergence, resulting in a combination which requires very few solver iterations. 
}


%\jnlcitation{\cname{%
%\author{Taylor M.},
%\author{Lauritzen P},
%\author{Erath C}, and
%\author{Mittal R}}.
%\ctitle{On simplifying ‘incremental remap’-based transport schemes.} \cjournal{\it J Comput Phys.} \cvol{2021;00(00):1--18}.}

\maketitle

\renewcommand\thefootnote{}
%\footnotetext{\textbf{Abbreviations:} ANA, anti-nuclear antibodies; APC, antigen-presenting cells; IRF, interferon regulatory factor.}

\renewcommand\thefootnote{\fnsymbol{footnote}}
\setcounter{footnote}{1}

\section{Motivation}

Linear systems (i.e., find $\M{x}$ solution of $\M{Ax = f}$ for a given square matrix $\M{A}$ and 
right-hand-side multi-column vectors $\M{f})$ are ubiquitous in science and engineering fields.
The development and utilization of highly efficient linear algebra software is crucial for solving linear systems in practical applications. There are two primary strategies to consider: iterative or direct solvers, each with its own advantages and drawbacks. Direct solvers aim for machine precision but can be memory-intensive, as the factorization process may generate significant fill-in without an effective reordering scheme. On the other hand, iterative solvers are often the preferred choice for very large sparse systems due to the low memory requirements of performing matrix-vector multiplication operations. 
However, their effectiveness relies heavily on preconditioners to improve convergence behavior.
%However, their effectiveness relies heavily on preconditioners, which improve the conditioning of the iterative process by providing a more easily applicable approximation of the original system.
%However, their effectiveness depends significantly on the use of preconditioners, which improve the conditioning of the iterative process using an approximation of the original system which is easier to apply. 
%involve solving linear systems that are simpler to manage than the original system. 
%However, their effectiveness depends significantly on the use of efficient preconditioners, which involve solving linear systems that are simpler to manage than the original system. 

While traditional solvers prioritize full precision, many applications do not require it.
 Moreover, with the decline of automatic performance improvements in single-threaded computing, parallel solvers have become indispensable for addressing large-scale problems \cite{Samehbook16}. 
Reducing precision can improve performance by lowering memory bandwidth demands and minimizing communication bottlenecks in distributed parallel environments. If higher accuracy is ultimately required, a reduced-precision solver can serve as a preconditioner for an iterative scheme, leading to what could be described as a hybrid or nearly direct solver. Certain applications, such as the FEAST eigenvalue solver \cite{PolizziFEAST,Polizzi2014,Gavin2018} or the TraceMIN eigenvalue solver \cite{Sameh:1982,Klinvex:2013}, already leverage reduced-precision computations where high accuracy for solving the resulting inner linear systems is unnecessary. 

Our work focuses on advancing the \Spike{}  algorithm \cite{Polizzi:2011}, a widely recognized parallel banded solver, by extending its capabilities to efficiently handle large, general sparse systems. We aim to develop a solver that is both highly parallel and adaptable to varying precision requirements. A central challenge is reducing communication overhead while preserving numerical accuracy, which is critical for solving large sparse linear systems effectively.
% going beyond its traditional application to banded systems. 
%We seek to develop a solver that is both highly parallel and flexible in accommodating diverse precision requirements.
%A key challenge in this endeavor is minimizing communication overhead without compromising numerical accuracy, which is %essential for the efficient solution of large sparse linear systems.
To address this problem, we propose a new family of \Spike{}  solvers and preconditioners, where the ``spikes'' generated during the \Spike{}  factorization process are computed using a low-rank approximation.
The resulting parallel, black-box \Spike{}  preconditioners for iterative solvers on distributed memory architectures not only significantly outperform block-Jacobi preconditioners but also provide robust support for diverse precision needs.

%A new family of SPIKE solver and preconditioner is proposed where 
%the "spikes" arising from the SPIKE algorithm 
%factorization stage are computed using a reduced-rank approximation.
%The resulting parallel black-box SPIKE preconditoners for iterative solvers on distributed memory computing architecture are not only significantly outperforming block-jacobi preconditioners, but are highly effective for accommodating diverse precision requirements.
%This work gives rise to a new family of SPIKE solvers for solving general sparse systems, 
%integrating reduced precision arithmetic to achieve an optimal trade-off between computational efficiency and accuracy. 

Following a background introduction to the \Spike{}  algorithm in Section~\ref{SPIKE_background}, Section 3 introduces a novel low-rank (LR) approximation approach tailored for \Spike{}. This innovation leads to the development of a new family of algorithms, termed \algname{}, which are detailed in Section 3.3. Section 4, dedicated to numerical experiments, encompasses three key components: first, a thorough qualitative case study; second, a quantitative performance analysis across various matrix systems; and third, an evaluation of parallel performance. This evaluation includes applying the new \algname{} schemes to system matrices generated by the NESSIE electronic structure code, as well as their integration with the FEAST solver.  

%    * Why parallel solvers?
%    
%       * "Free lunch" single threaded performance improvements are over
%       
%       * Distributed systems
%       
%    * Why reduced precision? 
%    
%        * Improve performance when full accuracy is not needed
%        
%        * Communication is often a bottleneck for distributed systems. Want to minimize communication requirements for a given level of accuracy. (Could probably cite Approximate Communication: Techniques for Reducing Communication Bottlenecks in Large-Scale Parallel Systems)
%        
%        * If greater accuracy is necessary in the end, may instead be used to precondition some %iterative scheme
        
%            * If the preconditioner is a good enough approximation, we might consider the overall scheme "hybrid" or "nearly direct."
            
%        * Some applications do not need high accuracy -- FEAST linear system solver
            
%            * I actually am not sure about the use of low precision linear solvers more generally.
    
\optpagebreak    
\section{\Spike{} background  }\label{SPIKE_background}

%\Spike{} is a parallel linear system solver designed for block tridiagonal matrices. 
%\todo{development history here}
%\todo{the intro to the previous paper would be perfect here...}
The \Spike{}  algorithm originates from the work of A. Sameh and D. Kuck in the late 1970s on tridiagonal systems \cite{Sameh:1978}, later extended to banded systems \cite{Chen:1978,Gallivan2012}. 
It is considered a domain decomposition method \cite{Eijkhout2012} for solving block tridiagonal systems. 
Unlike traditional \Mt{LU} factorization, \Spike{} introduces a \Mt{DS} factorization 
($\M{S}$ being the \Spike{}  matrix) that enhances parallel implementation by reducing communication costs. 
Over time, Sameh and his collaborators have developed various improvements and adaptations 
\cite{Dongarra:1984,Lawrie:1984,Berry:1988,Sarin:1999,Polizzi:2006,Polizzi:2007,Murat:2009,Murat:2010,Murat:2011,Spring2020} 
A key strength of \Spike{}  is its polyalgorithmic and hybrid nature, meaning it integrates multiple algorithmic approaches either direct or iterative, to optimize performance based on problem characteristics.  The algorithm using the $\M{DS}$ factorization achieves parallelism by decoupling the relatively large blocks along the diagonal of the $\M{D}$ matrix,  solving them independently, and then solving the spike matrix $\M{S}$ via the use of a smaller reduced system. 
Despite being logically divided into two stages (factorization and solve), computationally, the \Spike{}  algorithm comprises four main stages: (i) factorizing the diagonal blocks, (ii) computing the spikes,
(iii) solving the reduced system, (iv) retrieving the solution of the original system. 
Each of these stages can be accomplished in several ways, allowing a multitude of \Spike{}  schemes. 

Two notable variants are the recursive \Spike{}  algorithm for non-diagonally-dominant cases and the truncated \Spike{}  algorithm for diagonally-dominant cases. Depending on the variant, a system can be solved either exactly or approximately. When the system is solved exactly, \Spike{}  becomes a very effective dense banded solver that can significantly outperform the ScaLAPACK package on distributed memory architectures \cite{ScaLAPACK,Polizzi:2006}, as well as LAPACK on shared memory systems \cite{LAPACK,Spring2020}. The SPIKE-OpenMP package \cite{Mendiratta:2011,master,Spring2020,spike} is used as the default kernel for the FEAST eigensolver \cite{feast,FEASTUG} for solving eigenvalue problems that are banded.
GPU implementations of banded SPIKE have also been proposed by other authors \cite{Chang:2012,Venetis:2015}.

\Spike{} can also be used to solve general sparse systems, provided these systems are first reordered into a banded form using methods such as the reverse Cuthill–McKee algorithm \cite{Cuthill:1969} or weighted spectral ordering \cite{Murat:2009,Murat:2010}.
For many systems (e.g., those arising from 3D discretization schemes), the resulting band width may be very large while still maintaining sparsity within the band. 
%To reduce the computational cost of computing the full spikes (step (ii) above), two possible approaches have been employed, so far: (i) \Spike{} can function as either a dense or a more cost-effective sparse banded preconditioner when combined with iterative methods such as Krylov subspace techniques and iterative refinement \cite{Murat:2009,Murat:2010}; (ii) alternatively, the `\spike{}-on-the-fly' scheme \cite{Polizzi:2007} can be used for systems with very wide bands, where explicitly generating the spikes is difficult, but this increases the cost of solving the reduced system (step (iii) above). 
To reduce the computational cost of computing the full spikes (step (ii) above), two possible approaches have been employed, so far: (i) \Spike{} can function as a banded preconditioner when combined with iterative methods such as Krylov subspace techniques and iterative refinement \cite{Murat:2009,Murat:2010}; (ii) alternatively, the `\spike{}-on-the-fly' scheme \cite{Polizzi:2007} can be used for systems with very wide bands, where explicitly generating the spikes is difficult, but this increases the cost of solving the reduced system (step (iii) above). 
Having introduced the general background on the banded \Spike{} algorithm in this section, the next section will present a new \Spike{} framework designed for efficiently solving general sparse systems.

\subsection{SPIKE partitioning}
\Spike{} is a linear system solver, and so the problem to be solved is
\begin{equation}
\M{A}\M{x}=\M{f},
\end{equation}
where $\M{A}$ is an $n\times n$ matrix, $\M{f}$ is a collection of $n_{rhs}$ known vectors of size $n$ (the `right-hand-side'), and $\M{x}$ is a collection of $n_{rhs}$ unknown vectors.

\Spike{} is concerned with matrices which are block tridiagonal, and particularly effective for banded matrices. 
That is, matrices for which super and sub diagonals past some distance from the diagonal are equal to zero. 
Let us call that distance $\bw$ in both cases, resulting in a band width of $2\bw+1$. 
We do not lose generality by doing so; if the matrix has an unequal number of super and sub diagonals, we can simply imagine using the greater of the two to determine our $\bw$ value. 
In this case, a block tridiagonal partitioning may be imposed upon the matrix,
\begin{equation}
\M{A}\M{x} = 
\mybmat{cccc}{cccc}{
 \M{A}_{1}       & \M{\Bmat}_{1}        &                     &        \\
 \M{\Cmat}_{2}   & \M{A}_{2}              & \M{\Bmat}_{2}       &        \\
                 & \ddots               & \ddots              & \ddots \\
                 &                      & \M{\Cmat}_{p}       & \M{A}_{p}
}
\mybmat{c}{cccc}{
 \M{x}_1 \\
 \M{x}_2 \\
 \vdots \\
 \M{x}_p
}
=
\mybmat{c}{cccc}{
 \M{f}_1 \\
 \M{f}_2 \\
 \vdots \\
 \M{f}_p
}
\end{equation}
%\begin{equation}
%\M{A}\M{x} = 
%\left[\begin{array}{cccc}
%\M{A}_{1}&\M{\Bmat}_{1}&& \\
% \M{\Cmat}_{2}&\M{A}_{2}&\M{\Bmat}_{2}& \\
% &\ddots&\ddots&\ddots \\
% &&\M{\Cmat}_{p}&\M{A}_{p}
%\end{array}\right]
%\left[\begin{array}{cc}
% \M{x_1} \\
% \M{x_2} \\
% \vdots \\
% \M{x_p}
%\end{array}\right]
%=
%\left[\begin{array}{cc}
% \M{f_1} \\
% \M{f_2} \\
% \vdots \\
% \M{f_p}
%\end{array}\right]
%\end{equation}
%where a value of $p$ may be selected with the limitation that, for a matrix of size $n$, $p<n/\bw$, to retain the block tridiagonal structure.
where a value of $p$ may be selected more or less arbitrarily, with the limitation that too many partitions will result in a matrix which is no longer block tridiagonal.
Parallel performance will be extracted by associating each partition to a computational resource (MPI processes, in this implementation). 
The submatrices may be sized arbitrarily, with $\M{A}_i$ of size $n_i$. 
The only true restriction is $n_i > \bw$, although for load balancing purposes each partition should be roughly the same size. 
%and it is preferable to have $n_i \gg \bw$ for performance reasons.

Central to the \spike{} algorithm is the \Mt{DS} factorization.  
The matrix is factored into $\M{D}$, a block diagonal matrix, containing the submatrices $\M{A}_i$; and $\M{S}$, a block tridiagonal matrix. 
The partitioning is as follows,

%\footnote{Those familiar with the algorithm will notice that this matrix is typically called $\M{V_i}$. It was necessary to break with this convention, here, to prevent a name clash later when taking rank-reduced SVDs.} 

\begin{equation}
\M{A} = 
\M{DS} = 
%\left[\begin{array}{cccc}
\mybmat{cccc}{cccc}{
\M{A}_{1}&&& \\
 &\M{A}_{2}&& \\
 &&\ddots& \\
 &&&\M{A}_{p}
}
%\end{array}\right]
%\left[\begin{array}{cccc}
\mybmat{cccc}{cccc}{
\M{I}_{n_1}&\M{\Vmat}_{1}&& \\
 \M{\Wmat}_{2}&\M{I}_{n_2}&\M{\Vmat}_{2}& \\
 &\ddots&\ddots&\ddots \\
 &&\M{\Wmat}_{p}&\M{I}_{n_p}
}
%\end{array}\right]
,\label{DS_first}
\end{equation}
where $\M{I}_{n_i}$ denotes an identity matrix of size $n_i \times n_i$. 

The matrix $\M{D}$ will likely be familiar, as this is the matrix used by the well known Block Jacobi preconditioner \cite{lawn21}.
For a matrix with a reasonably small band width the matrix $\M{D}$ should contain the majority of the elements.
The matrix $\M{D}$ may be solved with perfect parallelism\textemdash there isn't any coupling between adjacent submatrices $\M{A_i}$, so each partition may be associated with a compute element and operations can be performed locally. 
$\M{S}$, on the other hand, could be viewed as representing necessary communication between neighboring partitions (in the next section, the use of low-rank approximations inside $\M{S}$ will be investigated, and the result is a preconditioner which bridges the gap between conventional \spike{} based solvers and the Block Jacobi preconditioner).  

%This results in a tw step solve operation for solving $\M{D}\M{S}\M{x}=\M{f}$:
The solve operation after \Mt{DS} factorization can be split into two stages,
\begin{align}
\M{D}\M{y}&=\M{f},\label{DSwith_solve1} \\  
\M{S}\M{x}&=\M{y}, \label{DSwith_solve2}  
\end{align}
where $\M{y}$ is an intermediary vector partitioned similarly to $\M{x}$. 

Due to the banded structure, the off-diagonal blocks of the original matrix, $\M{A}$, are of the form
\begin{equation}
\M{\Bmat}_{i} = 
%\left[\begin{array}{cc}
\mybmat{cc}{cc}{
 \M{0} & \M{0} \\
 \M{\bmat}_{i} & \M{0} \\
}
%\end{array}\right]
, 
\quad \text{ and } \quad
\M{\Cmat}_{i+1} = 
%\left[\begin{array}{cc}
\mybmat{cc}{cc}{
 \M{0} & \M{\cmat}_{i+1} \\
 \M{0} & \M{0} \\
}
%\end{array}\right]
, \text{ for } i\in 1...p-1,
\end{equation}
where $\M{\bmat}_{i}$ and $\M{\cmat}_{i}$ are $\bw \times \bw$ submatrices of $\M{\Bmat}_{i}$ and $\M{\Cmat}_{i}$ which contain all of their non-zero elements.
As a result note that the $\M{\Vmat}_i$ and $\M{\Wmat}_i$ submatrices in $\M{S}$ are mostly comprised of zeros, and the non-zero elements are restricted to tall, narrow ``spikes'' (from which the algorithm gets its name). 
\begin{align}
\M{\Vmat}_i = 
\inv{\M{A}}_i
%\left[\begin{array}{cc}
\mybmat{cc}{cc}{
 \M{0} & \M{0} \\
 \M{\bmat}_{i} & \M{0} 
} 
=
%\left[\begin{array}{cc}
\mybmat{cc}{c}{
 \M{\vmat}_i & \M{0} \\
}
%\end{array}\right]
,& \qquad
\M{\vmat}_i=
 \inv{\M{A}}_i
%\left[\begin{array}{c}
\mybmat{c}{cc}{
 \M{0}   \\
 \M{\bmat}_{i} 
}
%\end{array}\right]
, \qquad
\text{ for } i\in 1 \ldots p-1, \label{Veqn}
\\
\M{\Wmat}_{i} = 
 \inv{\M{A}}_{i}
%\left[\begin{array}{cc}
\mybmat{cc}{cc}{
 \M{0} & \M{\cmat}_{i} \\
 \M{0} & \M{0} 
} 
=
%\left[\begin{array}{cc}
\mybmat{cc}{c}{
 \M{0} & \M{\wmat}_{i} \\
}
%\end{array}\right]
,& \qquad
\M{\wmat}_{i}=
 \inv{\M{A}}_{i}
%\left[\begin{array}{cc}
\mybmat{c}{cc}{
 \M{\cmat}_{i} \\
 \M{0} 
}
%\end{array}\right]
, \qquad
\text{ for } i\in 2 \ldots p,
\end{align}
where $\M{\vmat}_{i}$ and $\M{\wmat}_{i}$ have the same number of columns, $\bw$, as $\M{\bmat}_{i}$ and $\M{\cmat}_{i}$.
These matrices, as well as the collections of vectors $\M{x}$ and $\M{y}$, may be broken up to extract the `tips' of each. 
\begin{equation}
\M{\wmat}_{i}=
%\left[
%\begin{array}{c}
\myaltbmatB{c}{ccc}{
 \M{\wmat}_{i,t} \\
 \M{\wmat}_{i,m} \\
 \M{\wmat}_{i,b} 
}
%\end{array}
%\right]
, 
\qquad
\M{\vmat}_{i}=
%\left[
%\begin{array}{c}
\myaltbmatB{c}{ccc}{
 \M{\vmat}_{i,t} \\
 \M{\vmat}_{i,m} \\
 \M{\vmat}_{i,b} 
}
%\end{array}
%\right]
, 
\qquad
\M{x}_i=
%\left[
%\begin{array}{c}
\myaltbmatB{c}{ccc}{
 \M{x}_{i,t} \\
 \M{x}_{i,m} \\
 \M{x}_{i,b} \\
}
%\end{array}
%\right]
, 
\qquad
\M{y_i}=
%\left[
%\begin{array}{c}
\myaltbmatB{c}{ccc}{
 \M{y}_{i,t} \\
 \M{y}_{i,m} \\
 \M{y}_{i,b} \\
}
%\end{array}
%\right]
, 
\end{equation}
where $\M{\vmat}_{i,t}$, $\M{\vmat}_{i,b}$, $\M{\wmat}_{i,t}$, $\M{\wmat}_{i,b}$, $\M{x}_{i,t}$, $\M{x}_{i,b}$, $\M{y}_{i,t}$ and $\M{y}_{i,b}$ all have $\bw{}$ rows.
This allows for the extraction of a reduced system. 
A three-partition example follows. 

\begin{equation}
\M{S}_{\red}\M{x}_{\red} 
=
\myaltbmatA{cccccc}{cccccc}{
\M{I}_{k}  &                   &\M{\vmat}_{1,t}&                  &               &                \\
           &    \M{I_{k}}      &\M{\vmat}_{1,b}&                  &               &                \\
           &    \M{\wmat}_{2,t}&\M{I}_k        &                  &\M{\vmat}_{2,t}&                \\
           &    \M{\wmat}_{2,b}&               &\M{I}_k           &\M{\vmat}_{2,b}&                \\
           &                   &               & \M{\wmat}_{3,t}  &\M{I}_k        &                \\
           &                   &               & \M{\wmat}_{3,b}  &               &\M{I}_{k}       \\
}
\myaltbmatA{c}{cccccc}{
 \M{x}_{1,t} \\
 \M{x}_{1,b} \\
 \M{x}_{2,t} \\
 \M{x}_{2,b} \\
 \M{x}_{3,t}\\
 \M{x}_{3,b}
}
=
\myaltbmatA{c}{cccccc}{
 \M{y}_{1,t} \\
 \M{y}_{1,b} \\
 \M{y}_{2,t} \\
 \M{y}_{2,b} \\
 \M{y}_{3,t}\\
 \M{y}_{3,b}
}
= \M{y}_{\red} \label{spikered} 
\end{equation}

%\end{array}\right]. 

This reduced system has a size of $(2\smul{}\bw{}\smul{}p) \times (2\smul{}\bw{}\smul{}p)$.
The fact that the reduced system size grows linearly with the number of partition interfaces could limit scalability. 
In addition, for matrices with very wide bands, forming the matrices $\M{\vmat}_{i}$ and $\M{\wmat}_{i}$ (required to explicitly extract the tips) may be impractical due to their size in memory. 
As a result, this is a frequent source of variant versions of the algorithm, which will be discussed in the next section. 

For now, let us assume that the reduced system may be solved somehow. 
Afterward, the overall solution may be recovered in parallel. 
Recovery may be performed either using the $\M{\vmat}_{i}$ and $\M{\wmat}_{i}$ matrices explicitly, or by using a solve operation, as seen in \eqns~(\ref{recovery1}-\ref{recovery2}). 

\begin{align}
 \M{x}_{i}
=
 \M{y}_{i} 
-
 \M{\wmat}_{i} 
\M{x}_{i-1,b}
-
 \M{\vmat}_{i} 
\M{x}_{i+1,t} 
=
 \inv{\M{A}_i}
\myaltbmatB{c}{ccc}{
\M{\cmat}_{i} \M{x}_{i-1,b} \\
\M{0}                     \\
\M{\bmat}_{i} \M{x}_{i+1,t}
}
\quad 
\text{ for } i\in 2 \ldots p-1.
\label{recovery1}
\end{align}

\begin{align}
 \M{x}_{1} 
=
 \M{y}_{1} 
-
 \M{\vmat}_{1} \M{x}_{2,t}
=
\M{y}_{1} 
-
 \inv{\M{A}_{1}}
\myaltbmatB{c}{cc}{
\M{0} \\ 
\M{\bmat}_{i} \M{x}_{2,t}
}
\M{x}_{2,t} \quad ; \quad
 \M{x}_{p} 
=
 \M{y}_{p} 
-
 \M{\wmat}_{p} 
\M{x}_{p-1,b} 
=
 \M{y}_{p} 
-
 \inv{\M{A}_{p}}
\myaltbmatB{c}{cc}{
\M{\cmat}_{i} \M{x}_{p-1,b} \\
\M{0}
}.
\label{recovery2}
\end{align}

Typically solve operations will be used for recovery. 
The matrices $\M{\vmat}_{i}$ and $\M{\wmat}_{i}$, are fairly large at $n_i \times \bw$ each, and dense, so their use may be inefficient (particularly for sparse matrices). 
In addition, using the solve-based form renders the submatrices $\M{\vmat}_{i,m}$ and $\M{\wmat}_{i,m}$, which contain most of the elements of $\M{\vmat}_{i}$ and $\M{\wmat}_{i}$, extraneous\textemdash the spikes must be formed in the factorization stage so their tips can be extracted, but subsequently the majority of their elements may be discarded, reducing overall memory usage.
Finally, in the next section we will see a technique to avoid explicitly forming the spikes at all. 

In summary, the steps of the solve operation for \spike{} are,
%\begin{enumerate}
\begin{align}
\text{Perform the D-stage:}      \spaceh  & \M{y}_i \leftarrow \inv{\M{A}_i} \M{f}_i                                         & \text{ for } & i\in 1...p  \label{D-stage-in-summary} \\ 
\text{Solve the reduced system:} \spaceh  & \M{x}_{\red} \leftarrow \inv{\M{S}}_{\red} \M{y}_{\red}                             &              &             \label{redsys-stage-in-summary} \\ 
\text{Recover solution:}         \spaceh  & \left\{ \begin{aligned} 
                                          \M{x}_1 &\leftarrow \M{y}_1                              - \M{\vmat}_{i}\M{x}_{2,t}    \\
                                          \M{x}_i &\leftarrow \M{y}_i - \M{\wmat}_{i}\M{x}_{i-1,b} - \M{\vmat}_{i}\M{x}_{i+1,t}  \\
                                          \M{x}_p &\leftarrow \M{y}_p - \M{\wmat}_{i}\M{x}_{p-1,b}                             
  \end{aligned}  \right. & \text{ for } & i\in 2...p-1 \label{recovery-stage-in-summary}
\end{align}
%\end{enumerate}

The D-stage and recovery may be performed with perfect parallelism, while the reduced system requires synchronization. 
In the next section, some strategies to ease the solution of the reduced system will be discussed. 

\subsection{\Spike{} reduced system solvers} \label{redsys_importance}
%-- Necessary for my part: \\
%-- -- \Spike OTF \\
%-- -- Iterative \\
%-- -- Recovery: \\ 
%-- -- -- Direct  \\
%-- -- -- Sweep   \\
%-- Generally interesting \\
%-- -- LU/UL           \\
%-- -- Special PARDISO factorization \\
%-- -- Recursive \Spike               \\ \\
%======================
%-- Outline, Method centric:\\
%-- -- Solving directly \\
%-- -- -- $2 \times 2$ \\
%-- -- -- Note about Recursive \Spike \\
%-- -- -- Truncated \\
%-- -- Solving iteratively \\
%-- -- -- \Spike OTF \\
%-- -- -- Note about LU/UL \\
%-- -- Recovery: \\ 
%-- -- -- Direct \\
%-- -- -- Sweep  \\
%\\
%======================
%-- Alternative Outline, decision centric. Property centric?: \\
%-- -- Solve the reduced system direcly or iteratively / Form V and W explicitly? \\
%-- -- -- Pros: \\ 
%-- -- -- -- Necessary if S to be solved directly \\
%-- -- -- -- Allows for extraction of tips \\
%-- -- -- Cons: \\
%-- -- -- -- V and W could be quite large \\
%-- -- -- -- Dense, even if V and W are sparse (note rearranging them) \\ 

%======================
% FREE-WRITING

The specifics of the operations used for \eqns~(\ref{D-stage-in-summary}) through~(\ref{recovery-stage-in-summary}) depend on the nature of the matrix being solved. 
An important aspect is the reduced system, because the size of the reduced system matrix grows linearly with the number of partitions, $p$ (specifically, the number of rows and the number of columns in the reduced system both grow linearly with $p$, but it is banded, so the growth of the number of elements is linear in $p$).
%To ameliorate this, variant algorithms have been developed for particular types of matrices. 
%Here, a quick overview of some of these variants will be provided.
To address the reduced system, various algorithmic variants have been developed for specific types of matrices. This section provides a brief overview of some of these approaches.

An appealing goal is to design a reduced system solver with parallel scalability. %, since the size of the reduced system grows with $p$. 
For example, if the reduced system solve operation achieves linear weak scalability, then the size of the reduced system is matched to the amount of parallelism, and so the reduced system can be viewed as a small (in the sense that $\bw$ is small) overhead, constant in $p$. 
Two important properties are: a matrix may be dense in the band or sparse, and it may or may not be diagonally dominant. 
For matrices which are diagonally dominant, it is possible to produce a block diagonal reduced system, which can be trivially solved in parallel. 
Matrices which are sparse banded will typically have a much wider band width than those which are dense in the band, due to the memory efficiency advantage of the former (which is to say, sparse storage allows for the consideration of matrices which would not fit in memory if stored in a dense format). 
Due to their wide bands, sparse matrices will be given special treatment.  

An important choice is whether the reduced system will be solved iteratively or directly. 
A direct solver provides ease of use to the end user\textemdash removing concerns such as convergence criteria, possible non-convergence, and unpredictable solve times.
On the other hand, iterative solvers naturally lend themselves to parallelism, as the matrix multiplication operation is relatively easy to parallelize, particularly in the case of banded matrices. 

In the case of two partitions, and therefore only one interface, the solution is trivial, 
\begin{align}
\mybmat{c}{cccc}{
 \M{x}_{1,t} \\
 \M{x}_{1,b} \\
 \M{x}_{2,t} \\
 \M{x}_{2,b} \\
}
\mybmat{cccc}{cccc}{
\M{I}_{k}  &                   &\M{\vmat}_{1,t}&        \\
           &    \M{I}_{k}      &\M{\vmat}_{1,b}&        \\
           &    \M{\wmat}_{2,t}&\M{I}_{k}      &        \\
           &    \M{\wmat}_{2,b}&               &\M{I}_{k}
}
=
\mybmat{c}{cccc}{
 \M{y}_{1,t} \\
 \M{y}_{1,b} \\
 \M{y}_{2,t} \\
 \M{y}_{2,b} \\
}.
\end{align}
It is sufficient to just solve the $2 \times 2$ block matrix in the center of this system, 
\begin{align}
\mybmat{c}{cc}{
 \M{x}_{1,b} \\
 \M{x}_{2,t} \\
}
\mybmat{cc}{cc}{
  \M{I}_{k}      &\M{\vmat}_{1,b} \\
  \M{\wmat}_{2,t}&\M{I}_{k}         \\
}
=
\mybmat{c}{cc}{
 \M{y}_{1,b} \\
 \M{y}_{2,t} \\
}. \label{2x2}
\end{align}
As can be seen in \eqn~(\ref{recovery2}), $\M{x}_{1,t}$ and $\M{x}_{2,b}$ can be obtained in the recovery stage. 
\Eqn~(\ref{2x2}) is the ``$2\times 2$ kernel'' \cite{Polizzi:2007,Spring2020}.
The $2\times 2$ kernel, restricted as it is to two partitions, can be through of as a primitive or building block for more complex schemes. 

In, \eqn~(\ref{spikered}) it can be observed that the blocks of $\M{S}_{\red}$ are connected via the submatrices $\M{\vmat}_{i,t}$ and $\M{\wmat}_{i,b}$.
For some matrices, these connecting submatrices may be ignored, leading to a variant known as ``truncated \spike{}''. 
In Ref.~\citenum{Mikkelsen:2009}, it is shown that for diagonally dominant matrices, the $\M{\vmat}$  and $\M{\wmat}$ spikes decay exponentially away from the diagonal, with a rate at least equal to the degree of diagonal dominance.
This behavior can also be observed for a certain class of Green's function matrices including quantum physics applications \cite{Lin2016}.
In these cases, $\M{\vmat}_{i,t}$ and $\M{\wmat}_{i,b}$ may be treated as zeros and ignored. 
The three-partition example is shown here,
\begin{equation}
\M{S}_{\red, trunc} 
=
\myaltbmatA{cccccc}{cccccc}{
\M{I}_{k}  &                   &               &                &                  &           \\
           &    \M{I}_{k}      &\M{\vmat}_{1,b}&                &                  &           \\
           &    \M{\wmat}_{2,t}&\M{I}_{k}      &                &                  &           \\
           &                   &               &\M{I}_{k}       &\M{\vmat}_{2,b}   &           \\
           &                   &               &\M{\wmat}_{3,t} &\M{I}_{k}         &           \\
           &                   &               &                &                  & \M{I}_{k} \\
}. \label{truncated_redsys}
\end{equation}
This results in a block diagonal matrix, and so the blocks may be solved independently and in parallel, achieving linear weak scalability. 

Other direct solvers for the reduced system exist. 
For example, recursive \spike{} observes that the reduced system is itself banded, and applies \spike{} to it, extracting some parallelism. 
This strategy has proven quite effect for matrices with narrow, dense bands, comparing favorably to state of the art solvers \cite{Polizzi:2006,Spring2020}.
However, for sparse matrices with very wide bands, just forming the spikes to extract the tips, required for direct solution, may be impractical. 

The alternative to solving the reduced system directly is to solve it iteratively. 
For example, a Krylov subspace algorithm may be used. 
In this case it is necessary to represent the action of multiplication by the reduced system. 
The reduced system is shown in \eqn~(\ref{spikered}).
On a per-parition basis, the multiplication operation is, 
\begin{align}
\mybmat{c}{cc}{
\M{y}_{1,t}  \\
\M{y}_{1,b}  
}
=  &
\mybmat{c}{cc}{
\M{x}_{1,t}  \\
\M{x}_{1,b}  
}
+
\mybmat{c}{cc}{
\M{\vmat}_{1,t}  \\
\M{\vmat}_{1,b}  
}
\M{x}_{2,t} 
\label{per-partition_redsys_first}
\\
\mybmat{c}{cc}{
\M{y}_{i,t}  \\
\M{y}_{i,b}  
}
=  &
\mybmat{c}{cc}{
\M{x}_{i,t}  \\
\M{x}_{i,b}  
}
+
\mybmat{c}{cc}{
\M{\wmat}_{i,t}  \\
\M{\wmat}_{i,b}  
}
\M{x}_{i-1,b}
+
\mybmat{c}{cc}{
\M{\vmat}_{i,t}  \\
\M{\vmat}_{i,b}  
}
\M{x}_{i+1,t} 
,\quad \text{ for } i\in 2 \ldots p-1.
\label{per-partition_redsys_mid}
\\
\mybmat{c}{cc}{
\M{y}_{p,t}  \\
\M{y}_{p,b}  
}
=  &
\mybmat{c}{cc}{
\M{x}_{p,t}  \\
\M{x}_{p,b}  
}
+
\mybmat{c}{cc}{
\M{\wmat}_{p,t}  \\
\M{\wmat}_{p,b}  
}
\M{x}_{p-1,b}
\label{per-partition_redsys_last}
\end{align}
To avoid explicitly creating these submatrices, the ``\spike{}-on-the-fly'' strategy has been developed \cite{Polizzi:2007}. 
For \spike{}-on-the-fly, the action of multiplication by the tips of the spikes is represented instead using sparse multiplication and solve operations. 
\begin{align}
\mybmat{c}{cc}{
\M{\vmat}_{i,t}  \\
\M{\vmat}_{i,b}  
}
\M{\tempA}_{i+1,t}
=&
\mybmat{ccc}{cc}{
\M{I}_{\bw} & \M{0} & \M{0}       \\
\M{0}       & \M{0} & \M{I}_{\bw} \\
}
 \inv{\M{A}}_{i}
 \mybmat{c}{cc}{
 \M{0}   \\
 \M{\bmat}_{i} 
}
 \M{\tempA}_{i+1,t}
 ,\spaceh  \text{ for } i\in 1...p-1 \label{OTF-1}
\\
\mybmat{c}{cc}{
\M{\wmat}_{i,t}  \\
\M{\wmat}_{i,b}  
}
\M{\tempA}_{i-1,b}
=&
\mybmat{ccc}{cc}{
\M{I}_{\bw} & \M{0} & \M{0}       \\
\M{0}       & \M{0} & \M{I}_{\bw} \\
}
\inv{\M{A}}_{i}
\underbrace{
\underbrace{
\mybmat{c}{cc}{
 \M{\cmat}_{i} \\
 \M{0}   
}
}_{n_{i} \times \bw}
\underbrace{\M{\tempA}_{i-1,b}}_{\bw \times n_{rhs}}}_{n_{i} \times n_{rhs}}
,\spaceh \text{ for } i\in 2...p \label{OTF-2}
\end{align}
The operations of~\eqns{}~(\ref{OTF-1}) and~(\ref{OTF-2}) may be performed from right to left.
This avoids the need to generate the $\M{\vmat}_{i}$ and $\M{\wmat}_{i}$ matrices explicitly inside the factorization stage. 
It does, however, require the use of a sparse solve operations inside the multiplication operation for the reduced system. 
A potential advantage comes from the number of vectors used for these solve operations\textemdash $\bw$ when solving for $\M{\vmat}_{i}$ or $\M{\wmat}_{i}$, but only $n_{rhs}$ when performing \spike{}-on-the-fly reduced system iterations. 
As such, this enables the use of \spike{} in cases where the submatrices $\M{\vmat}_{i}$ or $\M{\wmat}_{i}$ would not fit in memory. 
%As such, this option is most appealing if the matrix has a wide band relative to $n_{rhs}$.

The use of a sparse solve operation inside the reduced system multiplication may be costly, and so this strategy is most applicable to cases in which few iterations are necessary. 
For example, those in which the reduced system is well conditioned, or a loose convergence criteria is acceptable.
Note that, should it be possible explicitly form the spikes, the truncated spike reduced system can also be used to precondition these iterations. 

This concludes the description of the \spike{} variants which are relevant to this implementation. It should be noted that this is far from an exhaustive survey. % \todo{more sources}.
In the next section, we will look at approximating these operations to trade accuracy for performance.

\optpagebreak    
%\section{Building blocks for \algname{} * }
\section{Low Rank Approximation for \spike{}}

As described above, the method of solving the \spike{} reduced system should be tailored to the matrix being solved. 
For matrices with a narrow band width, the spikes ($\M{\vmat}_{i}$ and $\M{\wmat}_{i}$) may be explicitly formed and the required submatrices may be extracted, but when the matrix has a wide band, forming these matrices explicitly may be impossible due to memory constraints. 
We note that some sparse linear system solvers (PARDISO for example) are able 
to extract the relevant bottom and top spike tips without computing the whole spikes \cite{Murat:2009}. However, for extremely large sparse bands, even the tips become too large to generate.
%this becomes too challenging to perform these computations  on matrices with very large sparse bands.
%it is still not possible to explicitly extract the relevant submatrices of $\M{\vmat}_{i}$ and $\M{\wmat}_{i}$.
Instead one may use the ``\spike{}-on-the-fly'' technique, described in Section~\ref{redsys_importance}.
A drawback of this strategy is that it relies on sparse solve operations to represent the action of multiplication by the reduced system in an iterative solver, which can become computationally expensive if many iterations are needed.

In this section we will look at using a low rank singular value decompositions (low-rank SVDs) to produce an approximation of the 
whole spikes.

\subsection{SVD representation of the spikes} \label{SVD_rep_spike}

We would like to avoid the explicit construction of the large $\M{\vmat}_{i}$ and $\M{\wmat}_{i}$ matrices, while still producing representations amenable to submatrix extraction.
Later, in Section~\ref{Specific_schemes}, we will look at specific solver implementations.
In the meantime, it should be kept in mind that the approximate representations discussed here will be used to construct preconditioners in most cases, rather than direct solvers, and so we will accept some inaccuracy in the interest of improving performance. 

The matrices of interest are,

\begin{align}
  \M{{\vmat}}_{i} =
  \inv{\M{A}_{i}} 
    \mybmat{c}{ccc}{
		\M{0     } \\
		\M{0     } \\
		\M{\bmat}_{i,b} \\
    }
  \label{DiBaC}
  \spaceh \text{ and } \spaceh
  \M{{\wmat}}_{i} =
  \inv{\M{A}}_{i} 
    \mybmat{c}{ccc}{
		\M{\cmat}_{i,t} \\
		\M{0     } \\
		\M{0     } 
    }
   . 
\end{align}

As such, the low-rank SVD of each may be formed and partitioned as follows, 

\begin{align}
  \begin{split}
  \M{\inv{A_{i}}}   
    \mybmat{c}{ccc}{
		\M{0      } \\
		\M{0      } \\
		\M{\bmat}_{i,b} 
    }
  &=
    \mybmat{c}{ccc}{
		\M{\vmat}_{i,t} \\
		\M{\vmat}_{i,m} \\
		\M{\vmat}_{i,b} 
    }
  \approx 
  \M{\til{u     }}_{\M{\vmat}_{i}} 
  \M{\til{\Sigma}}_{\M{\vmat}_{i}}
  \M{\til{v     }}_{\M{\vmat}_{i}}
  =
    \mybmat{c}{ccc}{
		\M{\til{u}}_{\M{\vmat}_{i,t}} \\
		\M{\til{u}}_{\M{\vmat}_{i,m}} \\
		\M{\til{u}}_{\M{\vmat}_{i,b}} \\
    }
  \M{\til{\Sigma}}_{\M{\vmat}_{i}}
  \M{\til{v}}_{\M{\vmat}_{i}}
  \\
  &\approx
    \mybmat{c}{ccc}{
		\M{\til{u}}_{\M{\vmat}_{i,t}}\M{\til{\Sigma}}_{\M{\vmat}_{i}}\M{\til{v}}_{\M{\vmat}_{i}} \\
		\M{\til{u}}_{\M{\vmat}_{i,m}}\M{\til{\Sigma}}_{\M{\vmat}_{i}}\M{\til{v}}_{\M{\vmat}_{i}} \\
		\M{\til{u}}_{\M{\vmat}_{i,b}}\M{\til{\Sigma}}_{\M{\vmat}_{i}}\M{\til{v}}_{\M{\vmat}_{i}} 
    }
  =
    \mybmat{c}{ccc}{
		\M{\til{\vmat}}_{i,t} \\
		\M{\til{\vmat}}_{i,m} \\
		\M{\til{\vmat}}_{i,b} 
    }
    =
    \M{\til{\vmat}}_{i} 
  \end{split}
  \spaceh&
  \text{ for }i \in 2 \dots p, 
  \\
  \begin{split}
  \inv{\M{A}}_{i}   
      \mybmat{c}{ccc}{
		\M{   \cmat     }_{i,t} \\
		\M{    0        }       \\
		\M{    0        } 
      }
  &= 
      \mybmat{c}{ccc}{
		\M{      W}_{i,t} \\
		\M{      W}_{i,m} \\
		\M{      W}_{i,b} 
      }
  \approx
    \M{\til{u     }}_{\M{\wmat}_{i}}
    \M{\til{\Sigma}}_{\M{\wmat}_{i}}
    \M{\til{v     }}_{\M{\wmat}_{i}}
  = 
   \mybmat{c}{ccc}{
		\M{\til{u}}_{\M{\wmat_{i,t}}} \\
		\M{\til{u}}_{\M{\wmat_{i,m}}} \\
		\M{\til{u}}_{\M{\wmat_{i,b}}} 
   }
  \M{\til{\Sigma}}_{\M{\wmat}_{i}}\M{\til{v}}_{\M{\wmat}_{i}}
  \\ 
  &\approx 
   \mybmat{c}{ccc}{
		\M{\til{u}}_{\M{\wmat_{i,t}}}\M{\til{\Sigma}}_{\M{\wmat_{i}}}\M{\til{v}}_{\M{\wmat_{i}}} \\
		\M{\til{u}}_{\M{\wmat_{i,m}}}\M{\til{\Sigma}}_{\M{\wmat_{i}}}\M{\til{v}}_{\M{\wmat_{i}}} \\
		\M{\til{u}}_{\M{\wmat_{i,b}}}\M{\til{\Sigma}}_{\M{\wmat_{i}}}\M{\til{v}}_{\M{\wmat_{i}}} 
   }
  = 
   \mybmat{c}{ccc}{
		\M{\til{\wmat}}_{i,t} \\
		\M{\til{\wmat}}_{i,m} \\
		\M{\til{\wmat}}_{i,b} 
    }
  = 
  \M{\til{\wmat}}_{i} 
  \end{split}
  \spaceh&
  \text{ for }i \in 1 \dots p-1,
\end{align}
where the individual submatrices shown conform, e.g., 
$\M{\til{u}}_{\M{\wmat}_{i,t}}\M{\til{\Sigma}}_{\M{\wmat}_{i}}\M{\til{v}}_{\M{\wmat}_{i}} = \M{\til{\wmat}}_{i,t}$.
The number of singular values retained by each low-rank SVD will be called $n_{svd}$ (it will be assumed for convenience that all of the low-rank SVDs will retain the same number of singular values). 
As such, $\M{\til{u}}_{\M{\vmat}_{i}}$ and $\M{\til{u}}_{\M{\wmat}_{i}}$ will have dimensions of $n_{i} \times n_{svd}$, $\M{\til{\Sigma}}_{\M{\vmat}_{i}}$ and $\M{\til{\Sigma}}_{\M{\wmat}_{i}}$ are diagonal matrices with dimensions of $n_{svd} \times n_{svd}$, and $\M{\til{v}}_{\M{\wmat}_{i}}$ and $\M{\til{v}}_{\M{\wmat}_{i}}$ will have dimensions of $n_{svd} \times \bw$. 
The primary benefit is that we now have control over the number of columns in the large $n_i \times n_{svd}$ matrix, and may make tradeoffs between representation accuracy and size.

An appropriate value for $n_{svd}$ must be selected. 
While it is possible to bound the inaccuracy introduced by the truncation of the singular values (i.e., $|| \M{\vmat}_{i} - \M{\til{u}}_{\M{\vmat}_{i}} \M{\til{\Sigma}}_{\M{\vmat}_{i}} \M{\til{v}}_{\M{\vmat}_{i}} ||_2 \leq \Sigma_{{\M{\vmat}_{i}},(n_{svd}+1)}$ where $\Sigma_{{\M{\vmat}_{i}},(n_{svd}+1)}$ is the first excluded singular value), the impact after propagation through subsequent operations may not easily be predicted. 
As a result, the selection of an appropriate number of singular to retain should be made heuristically by the user. 
In particular, it may depend on the nature of the input matrix and desired trade-off between factorization time and accuracy (which may be heavily application dependent). 

When computing the SVD, there are two concerns. 
First, as noted previously, forming the matrices $\M{{\vmat}}_{i}$ and $\M{{\wmat}}_{i}$ may be impractical due to their size, which makes a matrix-free method appealing. 
However, the symbolic representation of these matrices, shown in \eqn~(\ref{DiBaC}), involves a sparse solve operation of size $n_i$.  
This may be computationally somewhat costly, and so it is preferable to avoid a method which might require many multiplications by $\M{{\vmat}}_{i}$ or $\M{{\wmat}}_{i}$. 
%This rules out many of the most popular methods of performing the SVD\todo{what are they again?}. 
%Instead, the randomized SVD algorithm of Halko, Martinsson, and Tropp~\cite{doi:10.1137/090771806} is used to compute the singular values and vectors.
Instead, the randomized SVD algorithm of Halko, Martinsson, and Tropp~\cite{doi:10.1137/090771806} is used.
%Specifically, because the singular spectra of the $\M{{\vmat}}_{i}$ and $\M{{\wmat}}_{i}$ have, for all of the matrices we have investigated so far, decayed quite quickly, is it possible to use the simple ``proto-algorithm'' of that work (and, if the singular values do not decay quickly, then the matrix is likely not a good candidate for this preconditioner, as the ability to select a small $n_{svd}$ is important to the overall performance of the solver). 
Specifically, our investigations have shown that the singular spectra of $\M{{\vmat}}_{i}$ and $\M{{\wmat}}_{i}$ tend to decay rapidly for all matrices examined so far. This suggests that the simple ``proto-algorithm'' from that work could be effectively applied.
%Conversely, if the singular values do not decay quickly, the matrix may not be well-suited for this preconditioner, as selecting a small  $n_{svd}$ is crucial for the solver's overall performance.
This implementation has been configured to perform only a small number of passes over the matrix (on the order of 2 to 4), and an over-sampling value of $\frac{1}{2}\times n_{svd}$ to $1\times n_{svd}$.

One might wonder why the singular values of $\M{{\vmat}}_{i}$ and $\M{{\wmat}}_{i}$ have been used. 
An appealing alternative could be to instead take the low-rank SVDs of $\M{{\bmat}}_{i}$ and $\M{{\cmat}}_{i}$, i.e., in the first case form $\M{\til{u}}_{\M{\bmat}_{i,b}}\M{\til{\Sigma}}_{\M{\bmat}_{i}}\M{\til{v}}_{\M{\bmat}_{i}} \approx \M{{\bmat}}_{i}$, and then use $\inv{\M{A}}_{i} \left[\M{0}, \M{\til{u}_{\bmat}}_{i}\right]^T$ to generate something analogous to $\M{\til{u}_{\vmat}}_{i}$. 
This would result in an easier singular value decomposition (as it would be performed on smaller matrices, and the solve operation would not be included inside). 
%However, by the Eckart-Young-Minsky theorem~\cite{doi:10.1093/QMATH/11.1.50}, the most accurate rank-reduced representation of $\M{{\vmat}}_{i}$ and $\M{{\wmat}}_{i}$, for a given rank, is the rank-reduced SVD of those matrices.  
However, by the Eckart-Young-Minsky theorem~\cite{doi:10.1093/QMATH/11.1.50}, the most accurate low-rank representation of $\M{{\vmat}}_{i}$ and $\M{{\wmat}}_{i}$, for a given rank, is the low-rank SVD of those matrices.  
Which is to say, constructing these matrices some other way will result in a worse approximation, or a larger value for $n_{svd}$.
%Either strategy could result in a viable solver, but because the factorization is reused for many solve operations in FEAST, we favor the highest quality, more costly factorization. 
An example using each strategy is presented and discussed in Section~\ref{Sherman5_subsec}.
%The degree to which the tradeoff between approximation efficiency and computational cost favors either strategy depends on the nature of the matrix being solved, as well as the performance needs of the overall application. 

This concludes the description of the SVD representation of the \spike{} matrices. 
%In the following two sections, we will look at two uses for this representation. 
The next two sections will explore two applications of this representation.
First, we will use it to perform the matrix multiplication operation for the reduced system, as a building block for an iterative scheme to solve this matrix.  
Second, we will use the SVD representation to construct a preconditioner for such an iterative scheme.
It may be desirable to mix-and-match these SVD based operations with the more conventional \spike{} operations. 
Some particular combinations will be discussed in Section~\ref{Specific_schemes}.

%\subsection{SVD Reduced System Multiplication}\label{redsysmul_sec}
\subsection{SVD Accelerated Reduced System Solve}\label{redsysmul_sec}

Performing the reduced system multiplication operation is very straightforward with the SVD representations of $\M{\vmat}_{i}$ and $\M{\wmat}_{i}$. 
Recall the 3-partition example of the reduced system is shown in \eqn~(\ref{spikered}).
Looking at the rows involved in the computation of $\left[\M{y}_{2,t},\M{y}_{2,b}\right]^T$, we may extract the operations which must be performed on a per-partition basis. 
In \Eqns~(\ref{per-partition_1}) and~(\ref{per-partition_2}), we see these operations, along with their SVD based approximations. 

\begin{align}
\mybmat{c}{cc}{
\M{y}_{i,t}  \\
\M{y}_{i,b}  
}
=  &
\mybmat{c}{cc}{
\M{x}_{i,t}  \\
\M{x}_{i,b}  
}
+
\mybmat{c}{cc}{
\M{\wmat}_{i,t}  \\
\M{\wmat}_{i,b}  
}
\M{x}_{i-1,b}
+
\mybmat{c}{cc}{
\M{\vmat}_{i,t}  \\
\M{\vmat}_{i,b}  
}
\M{x}_{i+1,t} 
\label{per-partition_1}
\\
\approx &
\mybmat{c}{cc}{
\M{x}_{i,t}  \\
\M{x}_{i,b}  
}
+
\mybmat{c}{cc}{
		\M{\til{u}}_{\M{\wmat}_{i,t}} \\
		\M{\til{u}}_{\M{\wmat}_{i,b}} 
}
	\M{\til{\Sigma}}_{\M{\wmat}_{i}}\M{\til{v}}_{\M{\wmat}_{i}}\M{x}_{i-1,b}
+
\mybmat{c}{cc}{
		\M{\til{u}}_{\M{\vmat}_{i,t}} \\
		\M{\til{u}}_{\M{\vmat}_{i,b}} 
}
	\M{\til{\Sigma}}_{\M{\vmat}_{i}}\M{\til{v}}_{\M{\vmat}_{i}}\M{x}_{i+1,t}
\label{per-partition_2}
\end{align}
In addition to the reduction in computation cost due to dimensionality reduction in \eqn~(\ref{per-partition_2}), this representation can also reduce communication cost.
The submatrices $\M{x}_{i+1,t}$ and $\M{x}_{i-1,b}$ are not available locally\textemdash they are local to neighboring nodes, and as such must be communicated, with a size of $\bw\times n_{rhs}$ each.
Instead, the SVD components $\M{\til{v}}_{\M{\wmat}_{i}}$ and $\M{\til{v}}_{\M{\vmat}_{i}}$ are supplied to the previous and next nodes respectively at factorization time. 
Performing these multiplications before communication reduces size of the transmissions from $\bw\times n_{rhs}$ to $n_{svd} \times n_{rhs}$ each.

%As mentioned previously, this multiplication operation will be used to approximate the reduced system, and the approximated reduced system will be solved using an iterative method.
%This multiplication operation will be used to approximate the reduced system, and the approximated reduced system will be solved using an iterative method.

%Next, we will look at a preconditioner, to accelerate the convergence of such a solver. 
%\subsection{Reduced System Preconditioner}\label{TruncTruncRedsys}

This matrix multiplication operation can be used to create a solver for the approximate reduced system, when paired with an appropriate iterative solver. 
Because the reduced system will be solved iteratively, it is helpful to have a preconditioner to accelerate convergence. 
In this case, we will use the truncated \spike{} variant as inspiration.
For truncated \spike{} the key observation is that, for some matrices, the elements of the submatrices $\M{\vmat}_{i}$ and $\M{\wmat}_{i}$ decay rapidly away from the diagonal.\cite{Mikkelsen:2009}.
As a result, the submatrices $\M{\vmat}_{i,t}$ and $\M{\wmat}_{i,b}$ may be treated as zero matrices in this case.
Removing these submatrices from the reduced system results in a block diagonal matrix, which may be solved in parallel.   
Let us call this new system $\M{\dtil{S}_{\red}}$.
A three-partition example is shown in \eqn~(\ref{red_sys_prec_mid_1}).

\begin{align}
  \M{\til{S}_{\red}} =&
	 \myaltbmatA{cccccc}{cccccc}{
		\M{I} &       & \M{\til{\vmat}}_{1,t} &             &             & \\
		&   \M{I}     & \M{\til{\vmat}}_{1,b} &             &             & \\
		& \M{\til{\wmat}}_{2,t} &    \M{I}    &             & \M{\til{\vmat}}_{2,t} & \\
		& \M{\til{\wmat}}_{2,b} &             &    \M{I}    & \M{\til{\vmat}}_{2,b} & \\
		&             &             & \M{\til{\wmat}}_{3,t} &   \M{I}     & \\
		&             &             & \M{\til{\wmat}}_{3,b} &             & \M{I}
    }
    \approx
      \myaltbmatA{cccccc}{cccccc}{
		\M{I} &       &                   &             &             & \\
		&   \M{I}     & \M{\til{\vmat}}_{1,b} &             &             & \\
		& \M{\til{\wmat}}_{2,t} &    \M{I}    &             &                   & \\
		&                   &             &    \M{I}    & \M{\til{\vmat}}_{2,b} & \\
		&             &             & \M{\til{\wmat}}_{3,t} &   \M{I}     & \\
		&             &             &                   &             & \M{I}
     }\label{red_sys_prec_mid_1}
  = \M{\dtil{S}_{\red}} 
	 \\
	 \inv{\M{\dtil{S}_{\red}}} = &
	\myaltbmatC{cccc}{cccc}{
		\M{I} &       &             &              \\
		& 
		\myaltbmatCII{cc}{cc}{
		  \M{I}     & \M{\til{\vmat}}_{1,b} \\
		  \M{\til{\wmat}}_{2,t} &    \M{I}   
		}^{-1}
		&             &              \\
		&             &    
		\myaltbmatCII{cc}{cc}{
		  \M{I}     & \M{\til{\vmat}}_{2,b} \\
		  \M{\til{\wmat}}_{3,t} &    \M{I}   
		}^{-1}
		&              \\
		&             &             &    \M{I}     
    } \label{red_sys_prec_mid_2}
\end{align}

There are two improvements to be made for this preconditioner. 
First, the block $2\times2$ submatrices along the diagonal of the matrix shown in \eqn~(\ref{red_sys_prec_mid_2}) may, for a matrix with a large band, be too large to practically form, factorize, and solve. 
This will be addressed by using a block inversion to exploit the structure of the matrices and the SVD representation of the off-diagonal blocks. 
Second, note that the components of these blocks are associated with different MPI processes, as indicated by their subscripts. 
So, applying these blocks these will require some communication, which we would like to minimize. 

The operations to be performed on a per-block basis are
\begin{align}
  \M{y}_{1,t}& = \M{g}_{1,t} \label{red_sys_1}
  \\
  \mybmatA{c}{cc}{
	 \M{y}_{i,b}        \\
	 \M{y}_{i+1,t} 
  }
& = 
\inv{
\mybmatA{cc}{cc}{
		  \M{I}     & \M{\til{\vmat}}_{i,b} \\
		  \M{\til{\wmat}}_{i+1,t} &    \M{I}   
}}
 \mybmatA{c}{cc}{
	 \M{g}_{i,b}        \\
	 \M{g}_{i+1,t}
  } 
  = 
  \mybmatA{cc}{cc}{
	 \M{I} + \M{\til{\vmat}}_{i,b}\inv{\M{M}_i}\M{\til{\wmat}}_{i+1,t} & (-\M{\til{\vmat}}_{i,b}\inv{\M{{M}_i}}) \\
	 (-\inv{\M{M}_i} \M{\til{\wmat}}_{i+1,t})                      &    \inv{\M{M}_i}   
   } 
  \mybmatA{c}{cc}{
	 \M{g}_{i,b}        \\
	 \M{g}_{i+1,t} 
   } \mbox{for } i \in 2\dots p-2 \label{red_sys_mid}\\
  \M{y}_{p,t}& = \M{g}_{p,t}, \label{red_sys_3}
\end{align}
where $\M{M}_i$ is the Schur complement in the bottom-right corner, $\M{M}_i = \M{I}-\M{\til{\wmat}}_{i+1,t}\M{\til{\vmat}}_{i,b}$.
$\M{M}_i$ may be rearranged for easier solution using the Woodbury matrix identity.

\begin{align}
  \inv{\M{M}_i} &= \underbrace{\inv{(\M{I}-\M{\til{\wmat}}_{i+1,t}\M{\til{\vmat}}_{i,b})}}_{\bw\times \bw}
  = (\inv{\M{I}+(-\M{\til{u}}_{\M{\wmat}_{i+1,t}}\M{\til{\Sigma}}_{\M{\wmat}_{i+1}})(\M{\til{v}}_{\M{\wmat}_{i+1}}\M{\til{u}}_{\M{\vmat}_{i,b}})\M{\til{\Sigma}}_{\M{\vmat}_{i}}\M{\til{v}}_{\M{\vmat}_{i}})} \\
  &= \M{I} + \M{\til{u}}_{\M{\wmat}_{i+1,t}}\M{\til{\Sigma}}_{\M{\wmat}_{i+1}}\underbrace{\inv{(\inv{(\M{\til{v}}_{\M{\wmat}_{i+1}}\M{\til{u}}_{\M{\vmat}_{i,b}})}-\M{\til{\Sigma}}_{\M{\vmat}_{i}}\M{\til{v}}_{\M{\vmat}_{i}}\M{\til{u}}_{\M{\wmat}_{i+1,t}}\M{\til{\Sigma}}_{\M{\wmat}_{i+1}})}}_{n_{svd}\times n_{svd}}\M{\til{\Sigma}}_{\M{\vmat}_{i}}\M{\til{v}}_{\M{\vmat}_{i}} \label{red_sys_prec_inv}
\end{align}

The matrices $\inv{(\M{\til{v}}_{\M{\wmat}_{i+1}}\M{\til{u}}_{\M{\vmat}_{i,b}})}$ and $\inv{(\inv{(\M{\til{v}}_{\M{\wmat}_{i+1}}\M{\til{u}}_{\M{\vmat}_{i,b}})}-\M{\til{\Sigma}}_{\M{\vmat}_{i}}\M{\til{v}}_{\M{\vmat}_{i}}\M{\til{u}}_{\M{\wmat}_{i+1,t}}\M{\til{\Sigma}}_{\M{\wmat}_{i+1}})}$ each have dimensions of $n_{svd} \times n_{svd}$.
Assuming a favorable SVD, these should be reasonably small, and so dense computations may be used.
The former must be explicitly inverted, while the latter may be formed and \Mt{LU} factorized for later solve operations. 
As a result, we have avoided inversion or solve operations involving large dense matrices.

The communication requirements of the preconditioner come from the off-diagonal submatrices in \eqn~(\ref{red_sys_mid}), which, when applied to the right hand side result in
$(-\inv{\M{M}_i} \M{\til{\wmat}}_{i+1,t}) \M{g}_{i,b}$ and $(-\M{\til{\vmat}}_{i,b}\inv{\M{M}_i})\M{g}_{i+1,t}$. 
To reduce the communication volume, $\M{g}_{i,b}$ and $\M{g}_{i+1,t}$ may be multiplied by $\M{\til{v}}_{\M{\wmat}_{i+1}}$ and $\M{\til{v}}_{\M{\vmat}_{i}}$, respectively, before sending.

The first is straightforward, the SVD representation may be used immediately,
\begin{align}
 \inv{\M{M}_i} \M{\til{\wmat}}_{i+1,t} \M{g}_{i,b} = 
 \inv{\M{M}_i} \M{\til{u}}_{\M{\wmat}_{i+1,t}}  \M{\til{\Sigma}}_{\M{\wmat}_{i+1}} (\M{\til{v}}_{\M{\wmat}_{i+1}} \M{g}_{i,b}),
\end{align}	
The second matrix requires a slight manipulation to expose $\M{\til{v}}_{\M{\vmat}_{i}}$ on the right of $\M{\til{\vmat}}_{i,b}\inv{\M{M}_i}$.
\begin{align}
 \M{\til{\vmat}}_{i,b}&\inv{\M{M}_i}\M{g}_{i+1,t} \\
& = \M{\til{u}}_{\M{\vmat}_{i,b}} \M{\til{\Sigma}}_{\M{\vmat}_{i}} \M{\til{v}}_{\M{\vmat}_{i}}(\M{I} + \M{\til{u}}_{\M{\wmat}_{i+1,t}}\M{\til{\Sigma}}_{\M{\wmat}_{i+1}}\inv{(\inv{(\M{\til{v}}_{\M{\wmat}_{i+1}}\M{\til{u}}_{\M{\vmat}_{i,b}})}-\M{\til{\Sigma}}_{\M{\vmat}_{i}}\M{\til{v}}_{\M{\vmat}_{i}}\M{\til{u}}_{\M{\wmat}_{i+1,t}}\M{\til{\Sigma}}_{\M{\wmat}_{i+1}})}\M{\til{\Sigma}}_{\M{\vmat}_{i}}\M{\til{v}}_{\M{\vmat}_{i}})\M{g}_{i+1,t}  \\
& = \M{\til{u}}_{\M{\vmat}_{i,b}} \M{\til{\Sigma}}_{\M{\vmat}_{i}} (\M{\til{v}}_{\M{\vmat}_{i}} + \M{\til{v}}_{\M{\vmat}_{i}}\M{\til{u}}_{\M{\wmat}_{i+1,t}}\M{\til{\Sigma}}_{\M{\wmat}_{i+1}}\inv{(\inv{(\M{\til{v}}_{\M{\wmat}_{i+1}}\M{\til{u}}_{\M{\vmat}_{i,b}})}-\M{\til{\Sigma}}_{\M{\vmat}_{i}}\M{\til{v}}_{\M{\vmat}_{i}}\M{\til{u}}_{\M{\wmat}_{i+1,t}}\M{\til{\Sigma}}_{\M{\wmat}_{i+1}})}\M{\til{\Sigma}}_{\M{\vmat}_{i}}\M{\til{v}}_{\M{\vmat}_{i}})\M{g}_{i+1,t} \\
& = \M{\til{u}}_{\M{\vmat}_{i,b}} \M{\til{\Sigma}}_{\M{\vmat}_{i}} (\M{I} + \M{\til{v}}_{\M{\vmat}_{i}}\M{\til{u}}_{\M{\wmat}_{i+1,t}}\M{\til{\Sigma}}_{\M{\wmat}_{i+1,t}}\inv{(\inv{(\M{\til{v}}_{\M{\wmat}_{i+1,t}}\M{\til{u}}_{\M{\vmat}_{i,b}})}-\M{\til{\Sigma}}_{\M{\vmat}_{i}}\M{\til{v}}_{\M{\vmat}_{i}}\M{\til{u}}_{\M{\wmat}_{i+1,t}}\M{\til{\Sigma}}_{\M{\wmat}_{i+1}})}\M{\til{\Sigma}}_{\M{\vmat}_{i}})\M{\til{v}}_{\M{\vmat}_{i}}\M{g}_{i+1,t} 
\end{align}	

%Consequently, a set of SVD based building blocks for \spike{} have been created. 
In conclusion, a set of SVD based building blocks for \spike{} have been created. 
For the most costly operations\textemdash matrix inversion and \Mt{LU} factorization\textemdash matrices larger than $n_{svd} \times n_{svd}$ have been avoided. 
Communication is reduced by pre-distributing the appropriate components of the low-rank SVD, essentially enabling the use of efficient lossy compression of the vectors to be communicated. 
Now we will use these building-blocks to construct actual preconditioners and solvers. 

\subsection{Selected \algname{} Algorithms} \label{Specific_schemes}

As described in Section~\ref{SPIKE_background}, there is a family of \spike{} solvers, mostly differentiated by the methods used to perform the following steps. 

\begin{align}
\text{Perform the D-stage:}      \spaceh  & \M{y}_i \leftarrow \inv{\M{A}_i} \M{f}_i                                         & \text{ for } & i\in 1...p  \label{D-stage-again} \\ 
\text{Solve the reduced system:} \spaceh  & \M{x}_{\red} \leftarrow \inv{\M{S}}_{\red} \M{y}_{\red}                             &              &          \label{redsys-stage-again1} \\ 
\text{Recover solution:}         \spaceh  & \left\{ \begin{aligned} 
                                          \M{x}_1 &\leftarrow \M{y}_1                              - \M{\vmat}_{i}\M{x}_{2,t}   \\
                                          \M{x}_i &\leftarrow \M{y}_i - \M{\wmat}_{i}\M{x}_{i-1,b} - \M{\vmat}_{i}\M{x}_{i+1,t} \\
                                          \M{x}_p &\leftarrow \M{y}_p - \M{\wmat}_{i}\M{x}_{p-1,b}                             
  \end{aligned}  \right. & \text{ for } & i\in 2...p-1 \label{recovery-stage-again2}
\end{align}

All our algorithms use Intel-MKL-Pardiso to solve the individual blocks, $\M{A_i}$, of Eq. (\ref{D-stage-again}) simultaneously.
With the set of low-rank SVD based building-blocks that we have constructed, many approximations of \eqns~\ref{redsys-stage-again1} and~\ref{recovery-stage-again2} would be possible. 
In table ~\ref{solver-list} some particularly interesting combinations are listed. 

%\newcommand{\newtilde}{\scalebox{.5}[.4]{\trimbox{0pt 1.8pt}{$\sim$}}}
%\newcommand{\newtildeb}{\scalebox{.5}[.4]{\trimbox{0pt 1.4pt}{$\sim$}}}
%\newcommand{\newtilde}{\scalebox{1.3}[1.3]{\trimbox{0pt 1.4ex}{\textasciitilde}}}
%\newcommand{\newtildeb}{\scalebox{1.3}[1.3]{\trimbox{0pt 1.3ex}{\textasciitilde}}}

% I had a little trouble picking the tilde symbol for the approximated matrices 
% in a way that was both pretty, portable, and required minimal packages. 
% In particular the double tilde is tricky.

%\newcommand{\til}[1]{\tilde{#1}{}}
%\newcommand{\til}[1]{\overset{ \scriptscriptstyle \sim}{#1}}
%\newcommand{\til}[1]{\accentset{ \scriptscriptstyle \sim}{#1}}
%\newcommand{\til}[1]{\stretchedtilde{#1}}
%\newcommand{\til}[1]{\stackengine{.5pt}{#1}{\scriptscriptstyle \sim}{O}{c}{F}{F}{S}}
%\newcommand{\til}[1]{\trlap[.5pt]{\scriptscriptstyle \sim}{#1}}
%\newcommand{\dtil}[1]{\tilde{\tilde{#1}}{}}
%\newcommand{\til}[1]{\stackon[.5pt]{#1}{\scriptscriptstyle \hspace{0pt} \sim}{}}
%\newcommand{\dtil}[1]{\stackon[-.5pt]{\stackon[.5pt]{#1}{\scriptscriptstyle \hspace{0pt} \sim}}{\scriptscriptstyle \hspace{0pt} \sim}{}}

%\begin{tabular}{|l|p{.2\linewidth}|p{.2\linewidth}|} \hline 
\begin{table}[htbp]
\renewcommand{\arraystretch}{1.3}
\centering 
\begin{tabular}{|c|c|c|} \hline 
{\bf Solver}                  & {\bf Reduced System}                                                                              & {\bf Recovery}        \\ \hline
\algname{}-I            & ${\M{\til{S}}_{\red}}$, solved iteratively, preconditioned with $\M{\dtil{S}}_{\red}$      & \tworow{Use $\M{\til{\vmat}}_{i}$ and/or $\M{\til{\wmat}}_{i}$ to approximate \eqn~\ref{recovery-stage-again2}}      \\ \cline{1-2} 
\algname{}-T            & $\M{\dtil{S}}_{\red}$, solved directly                                                      &  \\ \hline
\algname{}-OTF            & ${\M{{S}}_{\red}}$, solved iteratively, preconditioned using $\M{\dtil{S}}_{\red}$          & Performed exactly using sparse solve operations, as in \eqns~(\ref{recovery1}-\ref{recovery2}) \\ \hline
\end{tabular} 
\caption{Examples of \algname{} solvers. 
As a reminder, ${\M{{S}}_{\red}}$ is the true \spike{} reduced system (\eqn~\ref{spikered}), 
% $\M{\tilde{S}}_{\red}$ (\eqn~\ref{red_sys_prec_mid_1}), and $\M{\tilde{\tilde{S}}}_{\red}$ (\eqn~\ref{red_sys_prec_mid_2}). 
 $\M{\tilde{S}}_{\red}$ is the approximation of the reduced system using the low-rank SVD representation of the spikes (\eqn~\ref{red_sys_prec_mid_1}), and $\M{\tilde{\tilde{S}}}_{\red}$ is the further approximated version of the reduced system, based on truncated \spike{} (\eqn~\ref{red_sys_prec_mid_2}). We note that BiCGStab is used as the iterative solver. %%
}
\label{solver-list}
\end{table}

Because the low-rank approximations involve only small dense operations and involve very efficient communication operations, it is possible to borrow strategies from dense-in-the-band \spike{} solvers, exploiting the new efficient building blocks to apply these strategies to sparse-banded matrices.  
\algname{}-I and \algname{}-T represent this idea.
\algname{}-I is based on the original \spike{} algorithm, and solves the low-rank approximation of the true reduced system.
Combined with the use of $\M{\til{\vmat}}_{i}$ and $\M{\til{\wmat}}_{i}$ to recover the full solution, this preconditioner could be thought of as simply replacing the spikes in $\M{S}$ with the low-rank approximation. 
The level of approximation is controlled fully by the accuracy of the low-rank SVD.  
This preconditioner does introduce the issue of nested iterations\textemdash the reduced system must be solved to the working precision. 
However, in the experiments we have performed, inside convergence is usually reasonably fast (on the order of 10 or 20 inner iterations per outside iteration).  

If it is desirable to avoid the issue of nested iteration, \algname{}-T is an option.
In this case, the iterative solution of the reduced system is replaced by a direct use of the matrix $\M{\dtil{S}}_{\red}$. 
This produces a low-rank approximation of the ``truncated \spike{}'' solver. 
Truncated \spike{} is specialized to matrices for which the spikes decay away from the diagonal (diagonally dominant matrices being the most intuitive example). 
As such, depending on the nature of the matrix being solved, this may introduce an upper limit on the quality of this preconditioner\textemdash if the matrix is not a candidate for the conventional truncated spike solver, then our preconditioner can only represent it so accurately, whatever value is selected for $n_{svd}$.
However, this preconditioner is quite fast to apply and extremely scalable, and so it is often worthwhile to accept a slight loss in precontitioner quality (compared to \algname{}-I). 

Finally, we note that the \spike{} on-the-fly (\spike{}-OTF) algorithm, originally used  for cases in which the full-rank spikes are too large to generate explicitly, can also benefit from the new preconditioner $\M{\dtil{S}}_{\red}$. 
\Spike{} on-the-fly solves the true reduced system iteratively, with the symbolic representations of the spikes, shown in ~\eqns{}~(\ref{OTF-1}) and~(\ref{OTF-2}), used as a multiplication operation.
%Because \spike{}-OTF aims at addressing systems with very large sparse band widths, generating a preconditioner for this reduced system has remained a challenge until now.
Because \spike{}-OTF aims at addressing systems with very large sparse band widths, generating a preconditioner for this reduced system had been a challenge.
The use of the $\M{\dtil{S}}_{\red}$ preconditioner in \spike{}-OTF is a natural extension of this work and can significantly improve the convergence rate of the reduced system.
We will call this combination \algname{}-OTF.
It should be emphasized that by solving the true reduced system, \algname{}-OTF aims at producing a solution to the original problem (in contrast to the two other \algname{} schemes discussed above, which are primarily intended as preconditioners).
Consequently if the reduced system has converged to appropriate level of accuracy, outer iterations could be avoided entirely.

%previous attempts to generate a preconditoner for this reduced system have been <unsuccessful?>
%using only the top and bottom spikes,
%The use of the $\M{\dtil{S}}_{\red}$ preconditioner can indeed significantly improve the convergence rate of the reduced system. 
%an easy and straightforward improvement
%Once the reduced system has converged to appropriate level of accuracy, the full solution of the original linear system can be recovered without outer iterations.
%It can be difficult 

%The use of the $\M{\dtil{S}}_{\red}$ preconditioner for this reduced system is an easy and straightforward improvement, and can then can be solved
%with a higher level of accuracy.

%This is just a minor improvement to the existing \spike{} on-the-fly solver; the reduced system being solved by the iterative process is still the true \algname{} reduced system, and 

\section{Numerical Experiments}

Three sets of measurements are used investigate the effectiveness of this preconditioner. 
First, the impact of the preconditioner on a single matrix will be investigated in depth. 
Next, convergence behavior will be investigated for a collection of publicly available matrices \cite{Suitesparse}.
Finally, parallel performance on specific matrices arising from the electronic structure code NESSIE \cite{Kestyn:2020,NESSIE}, which inspired the development of this new family of \spike{} solvers, is reported. % the solver has been designed\textemdash matrices from the NESSIE \todo{..  description}. 

\subsection{Case Study: Sherman5} \label{Sherman5_subsec}

Sherman5 is a fairly old matrix, which was collected in the original Harwell-Boeing Sparse Matrix Collection \cite{10.1145/62038.62043}. 
This matrix was originally extracted from a finite difference model, used in an oil well simulation\cite{10.2118/13533-PA}.
It is reasonably well conditioned and quite small, so it is not very challenging to solve. 
However, it is indefinite, may be reordered to have a banded structure, and the small size makes it easy to work with, and so it is suitable to explore the effect of the preconditioner. 
For example, it is possible to explicitly create matrices such as $\M{\vmat}_{i}$ and $\M{\wmat}_{i}$, which must normally be avoided as they may be too large. 

To produce a banded structure, we employ the weight-aware reordering for band width reduction suggested in Ref.~\citenum{Murat:2010}. 
First, HSL-MC64 is used to produce column and row permutations which maximize the product of the elements along the diagonal, and to scale the diagonal to a magnitude of one (HSL-MC64 job option 5)\cite{HSL}. 
Next, the weighted spectral ordering (WSO) is applied, reducing the band width and pulling high-magnitude elements toward the diagonal. 
It should be noted that in this configuration HSL-MC64 is a non-symmetric reordering, and so it can remove zeros from the diagonal. 
In addition, the scaling operation may change the condition number of the matrix.

Finally, the process of reordering the matrix revealed some unconnected diagonal elements which do not contribute meaningfully to the solution. 
These lone elements are fully disconnected from the rest of the matrix (perhaps an artifact of unused boundary conditions in the application from which the matrix was originally taken, or something of that nature), and so they have been removed during reordering. 
After reordering, removal of unconnected elements, and scaling, the matrix has the following properties: 

\begin{itemize}
\item Condition number $3.2 \times \etothe{3}$ 
\item Number of rows and columns: $1638$
\item Band widths: $154$ (upper), $135$ (lower)
\item Percentage of weight on diagonal: $33\%$
\item Band density, $\#nonzeros/(band\ width\times n)$: $4\%$
\end{itemize}

\begin{figure}[h]
\centering 
\caption{Visualization of Sherman5 matrix before and after reordering. 
   Darker areas indicate a greater concentration of relatively large magnitude elements. A blur has been applied and magnitude is plotted on a log scale, to better show overall structure. }
\label{Sherman5_before_and_after}
  \subcaptionbox{Sherman5 with unconnected diagonal elements removed, but before further reordering.}%
  [.49\linewidth]{
  \includegraphics[keepaspectratio,width=\linewidth]{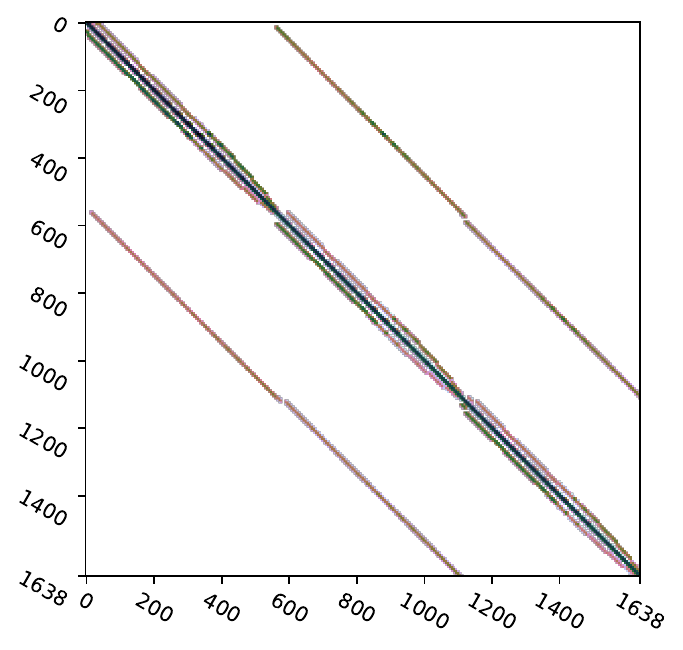}\hfill
  \label{Sherman5}}%
  \subcaptionbox{Sherman5 after HSL-MC64 and WSO.}%
  [.49\linewidth]{
  \includegraphics[keepaspectratio,width=\linewidth]{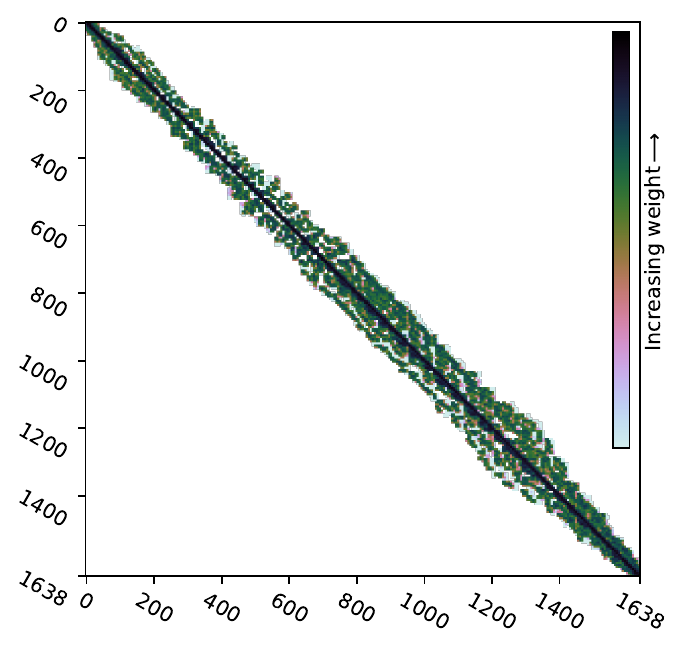}
  \label{Sherman5_wso}} 
\end{figure}

In \fig~\ref{Sherman5_before_and_after}, the matrix before and after reordering is shown. 
As can be seen, the band width has been reduced significantly, but not to the extent that a conversion to a dense in the band format could be seriously considered. 
\begin{figure}[htbp]
\centering
\caption{Characterizing the difference between using the SVD of the $\M{\vmat}_{i}$ and $\M{\wmat}_{i}$ versus $\M{\bmat}_{i}$ and $\M{\cmat}_{i}$. $\M{\vmat}_{i}$ and $\M{\wmat}_{i}$ result in a superior preconditioner due to faster singular value decay.}
\label{Sherman5_svd_B_vs_V}
  \subcaptionbox{Normalized singular values of the off diagonal submatrices for Sherman5, using three partitions. The singular values for $\M{\vmat}_{i}$ and $\M{\wmat}_{i}$, left, decay immediately in each case, while those of $\M{\bmat}_{i}$ and $\M{\cmat}_{i}$, right, often do not.\label{Sherman5_vwbcsvdComparison}}%
  [.49\linewidth]
  {
  \includegraphics[keepaspectratio,width=\linewidth]{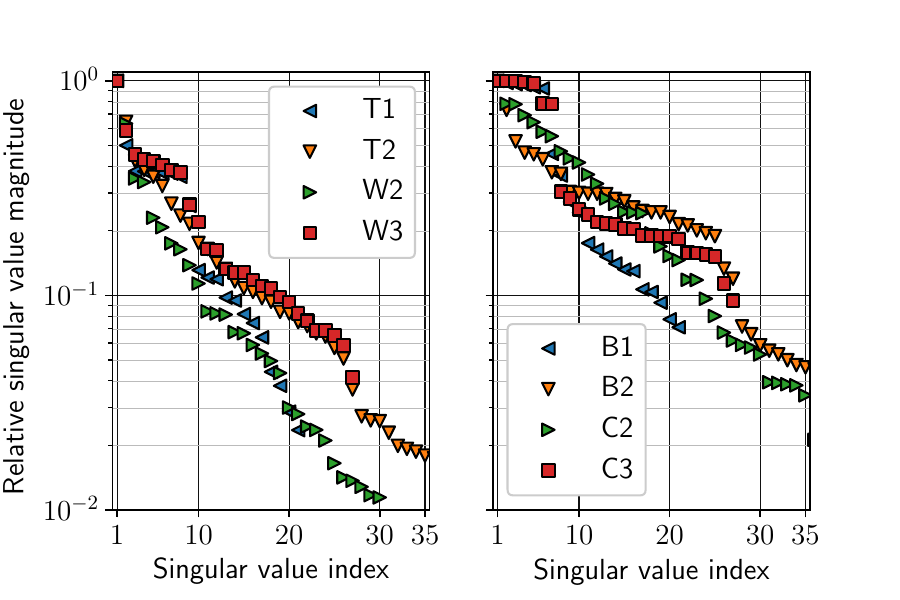}%
  }
  \hfill %\hspace{.05\linewidth}%
  \subcaptionbox{Comparing the quality of preconditioner produced by the two strategies using the condition number of the preconditioned matrix. Smaller numbers indicate better preconditioning. $n_{svd}=0$ corresponds to Block Jacobi. \label{Sherman5_preconditionedCond}}%
  [.49\linewidth]
  {
  \includegraphics[keepaspectratio,width=\linewidth]{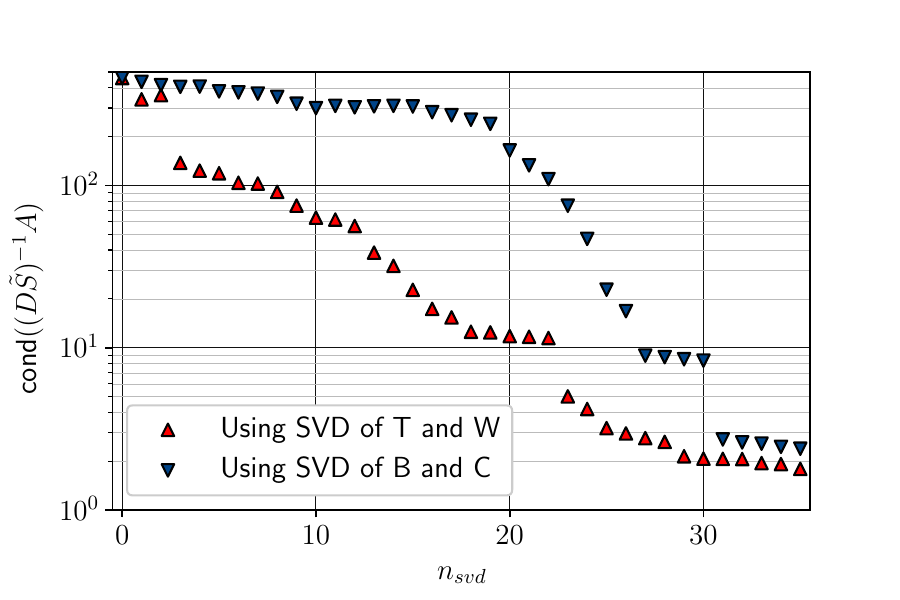}
  }%
\end{figure}
Let us now examine some qualities of the solver. 
The following measurements are described in terms of the quality achieved for a given $n_{svd}$.
Because many of the submatrices created scale in size with $n_{svd}$, minimizing this quantity will reduce computation requirements, as well as the communication requirements as seen in Section~\ref{redsysmul_sec}. 
In addition, a larger SVD size will increase the factorization time, as will be seen in Section~\ref{NessieMats}.

The Sherman5 matrix, permuted and scaled as described previously, has been broken up into three partitions and the \spike{} factorization has been performed. 
Three was selected for the number of partitions because any fewer would result in a special case\textemdash one partition is simply the base solver, while in the two partition case the $2\times2$ kernel can solve the reduced system directly. 
The first measurement shown relates to the quality of the approximation of the spikes produced by the low-rank SVD. 
In Section~\ref{SVD_rep_spike}, it is mentioned that the singular values of the matrices $\M{\vmat}_{i}$ and $\M{\wmat}_{i}$ decay fairly quickly.
An alternative strategy (mentioned briefly in Section~\ref{SVD_rep_spike}) to approximate $\M{\vmat}_{i}$ and $\M{\wmat}_{i}$ could be to use low-rank SVD representations of $\M{\bmat}_{i}$, and $\M{\cmat}_{i}$ as a starting point to generate the approximations of the spikes.
For example, $\M{{\bmat}}_{i} \approx \M{\til{u}}_{\M{\bmat}_{i}}\M{\til{\Sigma}}_{\M{\bmat}_{i}}\M{\til{v}}_{\M{\bmat}_{i}}$, and so $\M{\vmat}_{i} = \inv{\M{A}_{i}}[\M{0},\M{{\bmat}}_{i}^T]^T \approx (\inv{\M{A}_{i}}[\M{0}, \M{\til{u}}_{\M{\bmat}_{i}}^T]^T)\M{\til{\Sigma}}_{\M{\bmat}_{i}}\M{\til{v}}_{\M{\bmat}_{i}} $.
By the Eckart-Young-Minsky theorem, the approximation produced using the low-rank SVD of the actual spikes instead must be a better approximation than this, but computing the SVD of $\M{\bmat}_{i}$ and $\M{\cmat}_{i}$ would be easier, so it is worthwhile to check if the difference is significant. 

\Fig~\ref{Sherman5_vwbcsvdComparison} displays the relative magnitude of the singular values for matrices $\M{\vmat}_{i}$, $\M{\wmat}_{i}$, $\M{\bmat}_{i}$, and $\M{\cmat}_{i}$, for each partition.
The singular values have been normalized by the largest singular value in each set, to allow for easier comparison. 
Notice that the singular values of $\M{\vmat}_{i}$ and $\M{\wmat}_{i}$ drop in magnitude immediately, and quickly decay to around $10\%$ of their greatest value. 
The impact of this quick decay can be seen in \fig~\ref{Sherman5_preconditionedCond}, which shows the condition number of the matrix preconditioned by \algname{}-I, that is, $\inv{(\M{D} \M{\til{S}})}\M{A}$. 
The conditioning of this matrix immediately improves, corresponding to the drop in the singular values.  
As a result, it appears that there is a substantial advantage to be gained by approximating the actual spikes, rather than the matrices used to construct them.  

\begin{figure}[htbp]
\centering
\caption{Exploration of preconditioner quality for Sherman5 \phantom{Here is some phantom text, to deal with a bug in the Caption environment. It must be long enough to get us to two lines I guess?}}
 \label{Sherman5_svd_properties}
\subcaptionbox{Impact on preconditioner quality of varying both SVD size and number of partitions. \label{Sherman5_redsys_heatmaps}}%
 [.98\linewidth]
 {
 \includegraphics[keepaspectratio,width=\linewidth]{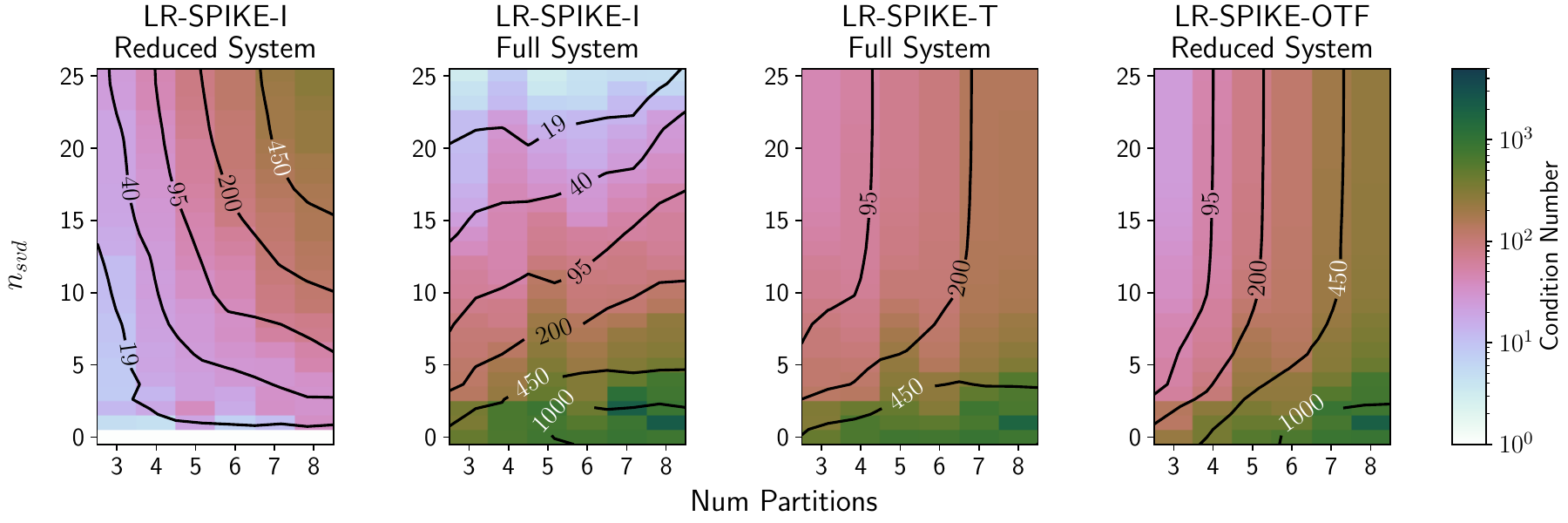}%
 }
 \subcaptionbox{Condition numbers associated with the three preconditioners for various SVD sizes.\label{Sherman5_redsys_cond}}%
 [.49\linewidth]
 {
 \includegraphics[keepaspectratio,width=\linewidth]{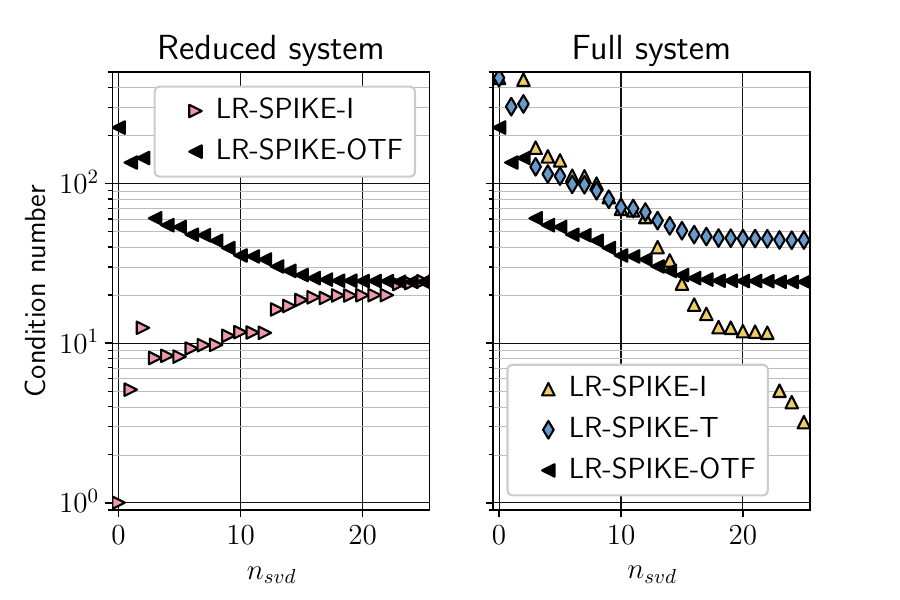}%
 }
 \hfill%\hspace{.05\linewidth}%
 \subcaptionbox{Effect of varying partition count on condition number, for a pair of SVD sizes. \label{Sherman5_redsys_cond_scaling}}%
 [.49\linewidth]
 {
 \includegraphics[keepaspectratio,width=\linewidth]{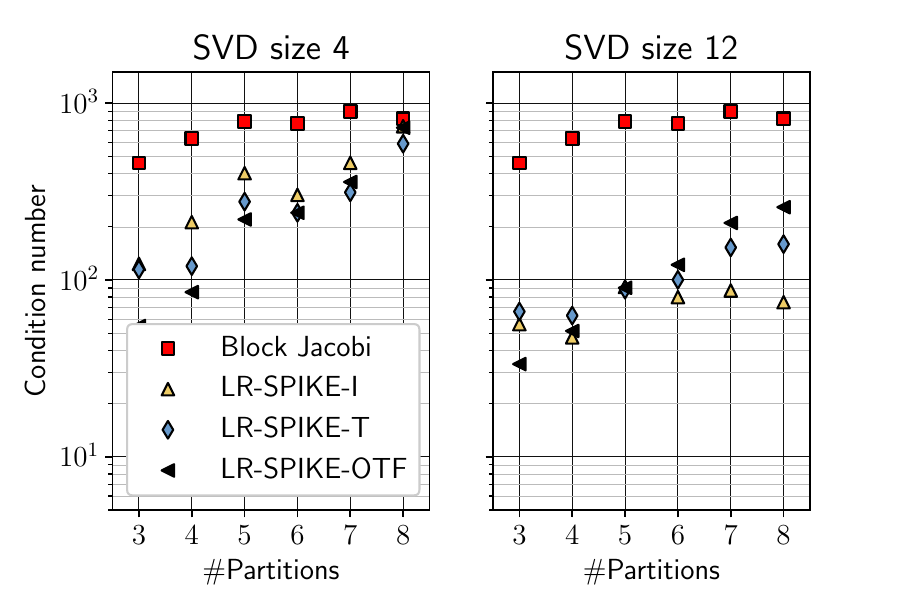}
 }%
% \hspace{.05\linewidth}%
\end{figure}

With the issue of which submatrices to approximate resolved, let us investigate the overall quality of preconditioners for the three example solvers discussed in Section~\ref{Specific_schemes}. 
The collection of measurements shown in \fig~\ref{Sherman5_svd_properties} explore the condition numbers associated with these solvers.
Larger condition numbers indicate that the iterative solver will have more difficulty solving the corresponding problem. 
The ideal condition number would be one, as this corresponds to an identity matrix; such a preconditioner could instead be used as a direct %%
solver. 

The quality of the preconditioner is primarily dependent on the number of singular values used to represent the spikes, and the number of partitions used. 
\Fig~\ref{Sherman5_redsys_heatmaps} shows a simultaneous exploration of these values, represented as heatmaps. 
Darker colors indicate higher condition numbers. 
Smoothed contour lines are also provided, to help identify overall trends.
The two leftmost heatmaps display condition numbers from \algname{}-I.
Uniquely among the three algorithms discussed here, \algname{}-I requires both inside and outside iterations; the approximate reduced system, which is solved iteratively, must be solved for each application of the full SPIKE preconditioner. 
Inside convergence is governed by the condition number of $\inv{\M{\dtil{S}}_{\red}}\M{\til{S}}_{\red}$, shown in the far left heatmap, outside by $\inv{\M{\til{S}}}\M{S}$, shown in the second heatmap from the left.
The next heatmap, second from the right, shows the preconditioner quality of the \algname{}-T preconditioner.
In this case the condition number of interest is that of $\inv{\M{\dtil{S}}}\M{S}$.
Measurements are only provided for the full problem, because the reduced system is solved directly.
Finally, for the \algname{}-OTF solver, the true reduced system is solved rather than an approximation. 
As such, outside iterations are not needed, and only the condition number of $\inv{\M{\dtil{S}}_{\red}}\M{S}_{\red}$ is shown. 

\Figs~\ref{Sherman5_redsys_cond} and~\ref{Sherman5_redsys_cond_scaling} show the impact of varying $n_{svd}$ and the number of partitions in isolation on these solvers. 
\Fig~\ref{Sherman5_redsys_cond} shows the condition numbers for a range of values of $n_{svd}$. 
The number of partitions has been held constant at three, so this could be considered a view of a vertical slice of the heatmaps. 
\Fig~\ref{Sherman5_redsys_cond}-left shows the condition numbers associated with the reduced systems, while \fig~\ref{Sherman5_redsys_cond}-right shows the effect of the preconditioner on the full system.  
Identical measurements are shown in both plots for \algname{}-OTF, again because this solver operates on the true reduced system, removing the need for outer iterations. 
For \algname{}-I, measurements of the condition number associated with the reduced system are present in \fig~\ref{Sherman5_redsys_cond}-left, while the condition number for the full system is shown in \fig~\ref{Sherman5_redsys_cond}-right.
Finally, as it was the case in \fig~\ref{Sherman5_redsys_heatmaps}, only the condition number associate with outside iterations is shown for the \algname{}-T preconditioner, because the approximation of the reduced system is solved directly in this case. 

In \fig~\ref{Sherman5_redsys_cond_scaling}, two values have been selected for $n_{svd}$, and the number of partitions has been varied.
These plots could  be considered a view of horizontal slices from the heatmaps. 
In addition to the three preconditioners discussed in this work, the conventional Block Jacobi preconditioner is included for comparison.
Because the preconditioner matrix of Block Jacobi is precisely the \spike{} matrix $\M{D}$, Block Jacobi could be seen as a limiting case of \spike{}, one in which $\M{\dtil{S}}$ and $\M{\til{S}}$ have been approximated by identity, in other words $n_{svd}=0$. 
For each of these preconditioners, increasing the number of partitions comes at the cost of decreasing preconditioner quality, because the number of elements in the approximated off-diagonal blocks increases.
In the case of Block Jacobi, as shown in \fig~\ref{Sherman5_redsys_cond_scaling}, the improvement over the unpreconditioned matrix (condition number $3.2\times \etothe{3}$) is fairly limited;  even for three partitions, a condition number of approximately $5\times \etothe{2}$ is shown, and the condition number levels off at just under $\etothe{3}$ when the number of partitions is large. 
For all of the \spike{} preconditioners, the quality for three partitions is much higher than that of Block Jacobi. 
The quality does still degrade as partitions are added, particularly in the cases of \algname{}-T and \algname{}-OTF. 
These two preconditioners are based on truncated \spike{}, which cannot fully represent the Sherman5 matrix (this matrix is not, for example, diagonally dominant, and so the spikes do not decay). 
As such, when $n_{svd}$ is sufficiently large, the nearness of the preconditioner to the original matrix is almost entirely determined by the degree to which truncated \spike{} can approximate the true matrix. 
In the case of the Sherman5 matrix, severely diminishing returns are seen after $n_{svd}=10$, and by $n_{svd}=16$ further increase becomes pointless. 
This can be seen most prominently in the two rightmost heatmaps of \fig~\ref{Sherman5_redsys_heatmaps}, as the contour lines become vertical in this range. 

For \algname{}-I, because the reduced system for the preconditioner is itself solved iteratively, both inner (reduced system) and outer (full system) convergence must be considered. 
The condition number associated with the reduced system for this preconditioner is that of $\inv{\M{\dtil{S}}_{\red}}\M{\til{S}}_{\red}$. 
As can be seen in \fig{}~\ref{Sherman5_redsys_cond}, this condition number begins at $1$ (for $n_{svd}=0$, both $\M{\dtil{S}}_{\red}$ and $\M{\til{S}}_{\red}$ are identity matrices) and increases with the number of singular values. 
With a large value selected for $n_{svd}$, the matrix $\M{\til{S}}_{\red}$ captures most of the true reduced system, and applying the preconditioner becomes nearly equivalent to solving the full problem. 
As a result, the condition number associated with the full system for \algname{}-I becomes quite small, with low single digits values for $n_{svd} > 22$.
Essentially, \algname{}-I provides the ability to arbitrarily increase the quality of the preconditioner for the full problem, at the cost of a more difficult iterative reduced system problem.
For the largest values of $n_{svd}$, the reduced system of \algname{}-I and the reduced system of \algname{}-OTF approach the same condition numbers. 

In summary, a detailed example of the effects of the two main tunable parameters\textemdash the number of singular values and the degree of parallelism\textemdash has been shown. 
Each of the new \algname{} preconditioners offers a significant improvement over the conventional Block Jacobi preconditioner. 
In general, a small value for $n_{svd}$ is favored, in the \algname{}-OTF and \algname{}-T cases because the underlying truncated spike approximation is only so accurate (outside of special cases such as diagonally dominant matrices), and in the \algname{}-I case to prevent the reduced system problem from becoming too difficult. 
\algname{}-OTF and \algname{}-T are particularly promising as purely direct preconditioners, while \algname{}-I can be used to provide additional accuracy if needed. 
Next, we will examine the convergence behavior of these solvers when used with a collection of SuiteSparse matrices. 

%
% Points Remaining: 
%  *  For -T and -OTF+T, the ultimate limit of the accuracy of the representation is a property of the matrix
%  *  Conventional truncated spike is strictly limited to handling matrices that it can represent well. Because we're operating as a preconditioner, we have more flexibility. 
%    *  Fully representing the truncated-spike based could be seen as a waste of precision. 
%  *  -I has mostly arrested the decline in preconditioner quality as the number of partitions grows
%    *  The spikes are represented reasonably well, and the cost has mostly been shifted to the reduced system 
%

\subsection{Preconditioner Effectiveness using SuiteSparse Matrices } 

To provide a general impression of the effectiveness of these new \spike{}-based preconditioners, they have been used to precondition BiCGStab iterations for a small collection of publicly available matrices. 
A summary of the test matrices is shown in \fig~\ref{matrix_table}.
The matrices selected are reasonably well conditioned, and small enough to be manipulated on a single machine. 
In addition, because this preconditioner is targeted towards general sparse banded matrices, symmetric positive definite and diagonally dominant matrices were intentionally avoided. 
Because \spike{} requires a matrix with a block tridiagonal structure, the reordering described in the previous section\textemdash HSL-MC64 followed by weighed spectral ordering\textemdash is used here as well.

The matrices tuma2 and bratu3D are from the SuiteSparse ``GHS\_indef'' collection \cite{Suitesparse}. 
These two matrices are symmetric, but indefinite. 
It should be noted that these are actually saddle point matrices, and as such it may be possible to tailor a custom preconditioner to exploit their particular structures\cite{benzi_golub_liesen_2005}. 
Here, we just use them as examples of indefinite matrices and perform the non-symmetric permutation mentioned above, destroying symmetry and structure, but generating banded matrices and avoiding problems associated with zeros on the diagonal.

The matrices waveguide3D, sme3D, and ns3Da are finite element matrices provided to Suitesparse from the COMSOL FEMLAB application. 
As indicated in Table~\ref{matrix_table}, before reordering these matrices already have most of their weight on the diagonal. 
They are general sparse matrices, but after reordering by WSO, they are suitably banded for use with these \spike{} preconditioners. 
An example, the reordering of sme3Da, is shown in \fig~\ref{sme3Da_before_and_after}.

\begin{figure}[htbp]
 \caption{Visualization of reordering for matrix sme3Da. \phantom{Here is some phantom text, to deal with a bug in the Caption environment. It must be long enough to get us to two lines I guess?}} \label{sme3Da_before_and_after}
%\captionsetup{width=0.8\linewidth,justification=raggedright}
\centering
  \subcaptionbox{Initial sme3Da matrix.}%
  [.49\linewidth]{\includegraphics[keepaspectratio,width=\linewidth]{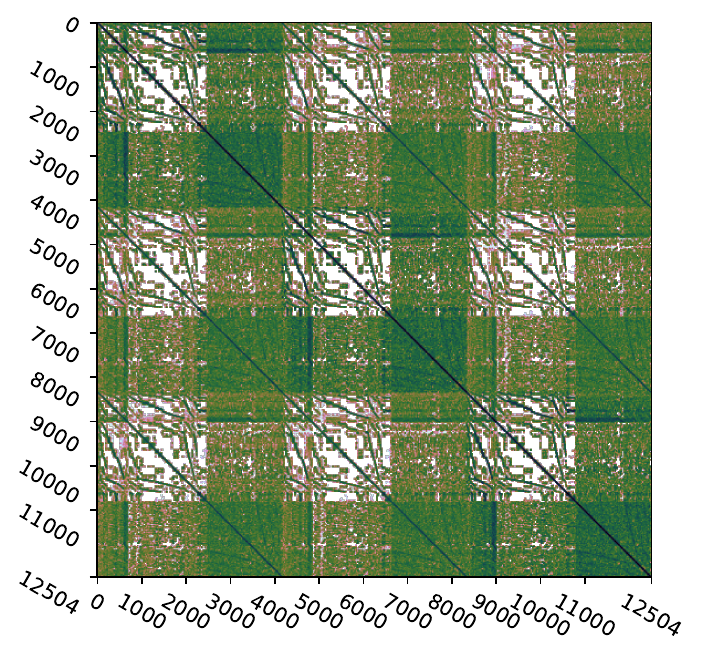}
  \label{sme3Da}}%
  \hfill %\hspace{.05\linewidth}%
  \subcaptionbox{sme3Da after HSL-MC64 and WSO are applied.}%
  [.49\linewidth]{\includegraphics[keepaspectratio,width=\linewidth]{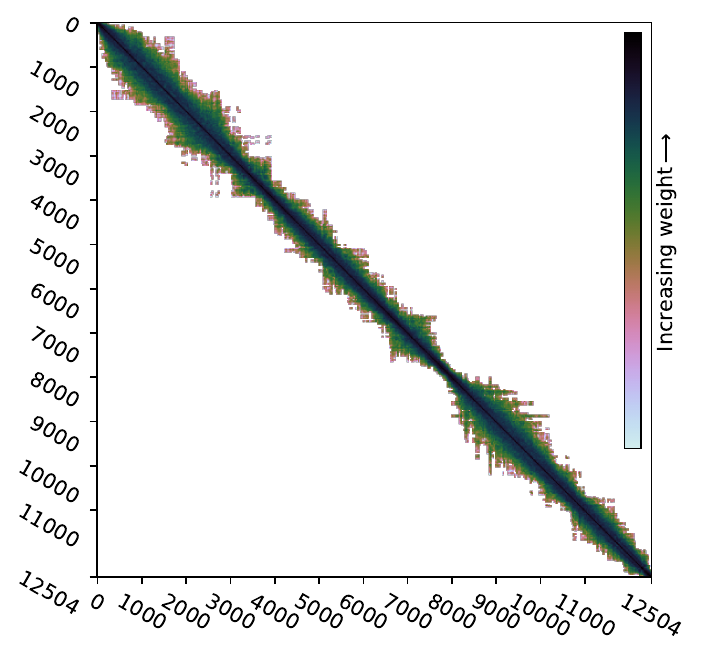}
  \label{sme3D_wso}}
\end{figure}

\begin{table}[htbp]
\centering 
\begin{tabular}{|ccc|ccc|ccc|} \cline{4-9}
  \multicolumn{3}{c|}{ }        & \multicolumn{3}{c|}{Before Reordering}       &   \multicolumn{3}{c|}{After Reordering}                                                \\        \hline
\tworow{Matrix} &\tworow{Size}     &\tworow{NNZ} & Condition                                     &   Symm-         & Diag.           & Condition                                     &  Half BW        & Diag.            \\        
                &                  &             &           Number                              &   etric?      &     Weight      &           Number                              &  (Up, Low)         &        Weight    \\ \hline 
sherman5        & 1638$^{\dagger}$ &  19119      &   $1.9 \times \etothe{5}$                     &    N           &     33\%        &   $3.2 \times \etothe{3}$                     &     (154,135)   &     33\%         \\ 
tuma2           & 12992            &  49365      &   $1.7 \times \etothe{3}$                     &    Y           &      4\%        &   $7.3 \times \etothe{3}$                     &  (1291,1051)    &     36\%         \\ 
bratu3D         & 27792            &  173796     &   $5.9 \times {\etothe{2}}                 $  &    Y           &     .2\%        &   $6.7 \times {\etothe{2}}                 $  &  (2211,2211)    &     49\%         \\ 
waveguide3D     & 21036            &  303468     &   $1.0 \times {\etothe{4}}                 $  &    N           &     32\%        &   $1.1 \times {\etothe{4}}                 $  &  (3466,3466)    &     32\%         \\ 
sme3Da          & 12504            &  874887     &   $5.2 \times {\etothe{7}}                 $  &    N           &     21\%        &   $1.1 \times {\etothe{7}}                 $  &  (1366,1366)    &     21\%         \\ 
ns3Da           & 20414            & 1679599     &   $7.1 \times {\etothe{2}}                 $  &    N           &     14\%        &   $7.1 \times {\etothe{2}}                 $  &  (3057,3057)    &     14\%         \\ 
\hline
%\multicolumn{6}{l}{\footnotesize{$\dagger$ disconnected diagonal elements were removed during reordering}} 
\end{tabular} 
\caption{Test matrix properties before and after reordering. $\dagger{}$ indicates fully disconnected diagonal elements were removed during the reordering process. \phantom{Here is some phantom text, to deal with a bug in the Caption environment. It must be long enough to get us to two lines I guess?}} \label{matrix_table}
\end{table}

Table~\ref{SuiteSparseBiCGStab} lists the number of BiCGStab iterations required to reach convergence while using each preconditioner. 
BiCGStab has been configured with a $\emin{7}$ convergence threshold.
An iteration limit of $100,000$ has been imposed, although it was only reached on one case.
For the ILU preconditioned solvers, as well as the unpreconditioned solver, the GNU Octave environment is used\cite{Octave}. 
Due to their distributed nature, the \algname{} preconditioners and Block Jacobi have been paired with an in-house distributed memory BiCGStab implementation. 

\begin{table}[htbp]
\centering 
\begin{tabular}{|cc|cc|cc|cc|cc|cc|cc|} \hline
\multicolumn{2}{|c|}{Preconditioner}  &   \multicolumn{12}{c|}{Matrix}                                                                                                                                         \\ \hline     
  Type               & $p_{q}$     &        \mco{sherman5}     & \mco{tuma2}        &     \mco{ bratu3D}          & \mco{waveguide3D}                      & \mco{sme3Da}           & \mco{ns3Da}              \\ \hline 
% \Spike{}-OTF        & \tworow{N/A}& \mcttworow{N/A}{22.5}  &\mcttworow{N/A}{159.5} &     \mcttworow{N/A}{19.5}   &     \mcttworow{N/A}{301}               &     \mcttworow{N/A}{795}& \mcttworow{N/A}{35}     \\  \cdashline{1-14}
 \Spike{}-OTF        &         N/A & \mct{N/A}{22.5}           & \mct{N/A}{159.5}   &     \mct{N/A}{19.5}         &     \mct{N/A}{301}                     &     \mct{N/A}{795}     & \mct{N/A}{35}     \\  \cdashline{1-14}
%                     &  16         &    \mct{N/A}{6.5}         & \mct{N/A}{103.5}   &     \mct{N/A}{14.5}         &     \mct{N/A}{192}                     &     \mct{N/A}{67.5}    &     \mct{N/A}{27.5}      \\ 
%\tworow{\Spike{}-OTF}&  32         &    \mct{N/A}{5 }          & \mct{N/A}{80.5}    &     \mct{N/A}{11.5}         &     \mct{N/A}{323.5}                   &     \mct{N/A}{101.5}   &     \mct{N/A}{25}        \\        
%                     &  64         &    \mct{N/A}{5 }          & \mct{N/A}{47}      &     \mct{N/A}{9.5}          &     \mct{N/A}{82}                      &     \mct{N/A}{61.5}    &     \mct{N/A}{17.5}      \\         
%                     &  128        &    \mct{N/A}{13}          & \mct{N/A}{11}      &     \mct{N/A}{8.5}          &     \mct{N/A}{7 }                      &     \mct{N/A}{19}      &     \mct{N/A}{14}        \\ \hline
                     &  16         &    \mct{N/A}{6.5}         &  \mct{N/A}{135}    &     \mct{N/A}{16}           &     \mct{N/A}{203}                     &     \mct{N/A}{75}      &     \mct{N/A}{28.5}      \\
\tworow{\algname{}-OTF}&  32         &    \mct{N/A}{5}           &  \mct{N/A}{104}    &     \mct{N/A}{11.5}         &     \mct{N/A}{319.5}                   &     \mct{N/A}{90.5}    &     \mct{N/A}{25}        \\
                     &  64         &    \mct{N/A}{5}           &  \mct{N/A}{51.5}   &     \mct{N/A}{9.5}          &     \mct{N/A}{95}                      &     \mct{N/A}{67}      &     \mct{N/A}{17.5}      \\
                     &  128        &    \mct{N/A}{9.5}         &  \mct{N/A}{12}     &     \mct{N/A}{8.5}          &     \mct{N/A}{7}                       &     \mct{N/A}{23}      &     \mct{N/A}{13.5}      \\ \hline
                     &  16         &    \mct{5}{64.5}          & \mct{52.5}{825.5}  &     \mct{13}{207}           &     \mct{204.5}{2417}                  &     \mct{66.5}{886}    &     \mct{28}{522}        \\
 \tworow{\algname{}-I}&  32         &    \mct{2}{26}            & \mct{48.5}{728.5}  &     \mct{10}{166}           &     \mct{285.5}{3649}                  &     \mct{62.5}{926}    &     \mct{23.5}{434.5}    \\
                     &  64         &    \mct{1.5}{19.5}        & \mct{22}{439}      &     \mct{8.5}{139}          &     \mct{73.5}{889.5}                  &     \mct{33}{615.5}    &     \mct{15.5}{278}      \\
                     &  128        &    \mct{0.5}{13.5}        & \mct{2}{31}        &     \mct{6.5}{105.5}        &     \mct{4}{46}                        &     \mct{8.5}{275.5}   &     \mct{12}{218}        \\ \hline
                     &  16         &    \mco{6}                & \mco{54}           &     \mco{14}                &     \mco{210.5}                        &     \mco{72.5}         &     \mco{26.5}           \\
 \tworow{\algname{}-T}&  32         &    \mco{4.5}              & \mco{55.5}         &     \mco{10}                &     \mco{300.5}                        &     \mco{55.5}         &     \mco{21.5}           \\
                     &  64         &    \mco{4.5}              & \mco{25.5}         &     \mco{8}                 &     \mco{76.5}                         &     \mco{33.5}         &     \mco{17}             \\
                     &  128        &    \mco{6.5}              & \mco{5}            &     \mco{7.5}               &     \mco{6}                            &     \mco{15.5}         &     \mco{14}             \\ \hline
Block                & \tworow{N/A}& \mco{\tworow{21.5}}   & \mco{\tworow{95.5}}    &    \mco{\tworow{22.5}} &     \mco{\tworow{160.5}} &   \mco{\tworow{168.5}} &     \mco{\tworow{32}}    \\        
   Jacobi            &             &&                       &&                         &&                       &&                         &&                       &&                         \\ \hline 
ILU$_0$              & N/A         &  \mco{ 22.5}        & \mco{ SF     }            &\mco{  23}              & \mco{3055.5  }           & \mco{SF      }         & \mco{FF      }           \\ \hline     
                     & $\emin{1}$  &    \mco{ 31.5}        & \mco{     267}          &\mco{  21}              & \mco{64409   }           & \mco{SF      }         & \mco{  3161.5}           \\             
ILU$_t$              & $\emin{2}$  &    \mco{ 11.5}        & \mco{       9}          &\mco{13.5}              & \mco{NC      }           & \mco{     469}         & \mco{    13.5}           \\             
                     & $\emin{3}$  &    \mco{  5.5}        & \mco{     5.5}          &\mco{ 7.5}              & \mco{    99.5}           & \mco{      66}         & \mco{     7.5}           \\ \hline             
None                 & N/A         &  \mco{251.5}        & \mco{ SF     }            &\mco{54.5}              & \mco{20212.5 }           & \mco{SF      }         & \mco{SF      }           \\ \hline     
\end{tabular}
\caption{BiCGStab iterations required to reach $\emin{7}$ residual norm. 
For \algname{}-I and \algname{}-T the value of $n_{svd}$ is listed in column $p_{q}$. All of the \spike{}-based preconditioners, as well as Block Jacobi, have been configured to use three partitions. 
Two ILU variants were used, ILU$_0$ (incomplete LU with no fill in), and ILU$_t$ (ILU with threshold based dropping of elements, using diagonal boosting to handle zero-pivots). 
For ILU$_t$, the dropping threshold is listed in column $p_{q}$. 
BiCGStab iterations involve two applications of the preconditioner and two convergence checks, so half-iterations are shown in some cases. 
For \algname{}-I, inside BiCGStab iterations (shown in parentheses) are required to solve the reduced system. 
The values shown represent the total count of inside BiCGStab iterations used in solving the full problem, aggregated across outside iterations.
A BiCGStab tolerance of $\emin{16}$ was selected for the \algname{}-I reduced system.
For the 'on-the-fly' variants, the reduced system BiCGStab solver targets a residual norm of $\emin{8}$ by default. 
If this does not result in sufficient accuracy for the overall solution, the reduced system residual norm tolerance is dynamically lowered and BiCGStab is restarted, using the previous reduced system solution as an initial guess. 
FF=Factorization Failure. 
SF=Failure in the solve stage (BCGStab stagnation detected, for example). 
NC=Did not converge in 100,000 steps.} 
\label{SuiteSparseBiCGStab}
\end{table}

For the \spike{}-based preconditioners increasing $n_{svd}$ improves the preconditioner quality in most cases. 
A very specific problem occurs with the matrix Sherman5. 
In this case, the value of $n_{svd}$ is close to the number of off-diagonals of the matrix, and a large number of near zero values appear in the singular value decompositions.
This tends to produce a poorly conditioned $2 \times 2$ reduced system preconditioner. 
Practically, it is easy to detect that the configured value of $n_{svd}$ is too large by checking for and automatically reduce it. 

In addition, all of the \spike{}-based solvers appear to have trouble with waveguide3d specifically at $n_{svd}=32$.
This also is the result of a poorly conditioned $2 \times 2$ reduced system preconditioner, although there is no obvious cause in this case. 
However, the challenge is overcome as $n_{svd}$ is increased.

Comparing \algname{}-I and \algname{}-T, the former reduces the number of outside iterations more effectively than the latter in almost every case, at the cost of inner iterations. 
The number of inner iterations for \algname{}-I is quite large in some cases, but it should be kept in mind that the computational and communication costs of these inner iterations scale linearly with $n_{svd}$. 
Note that the sizes selected for $n_{svd}$ are quite small relative to the half band widths of the matrices, less than 10\% for all matrices other than Sherman5. 

For \spike{}-OTF and \algname{}-OTF, only reduced system iterations are present. 
Because each reduced system iteration requires a large solve operation, these iterations should be considered as roughly comparable in cost to an outside iteration for \algname{}-I.
And in fact, the iteration counts observed for \algname{}-OTF are very similar to those of \algname{}-T in many cases, agreeing with the trends seen in Section~\ref{Sherman5_subsec}, \fig{}~\ref{Sherman5_svd_properties}. 

In almost every case, the \spike{}-based methods produce more effective preconditioners than the popular Block Jacobi and ILU$_0$ schemes.
Typically, they are similar in quality to ILU$_{t}$ with a drop threshold of $\emin{2}$ or $\emin{3}$, indicating that high quality preconditioning may be achieved with this technique.

In the next section, we will turn our attention to parallel scalability performance and timing results using a pair of problems extracted from our NESSIE application. 
%These matrices have a strong diagonal tendency, and so we will focus on  the \algname{}-T preconditioner. \todo{I think "strong diagonal tendency" is not quite the right description but I should just align this with the NESSIE matrix description later.}

\subsection{Parallel Performance}\label{NessieMats}

From the NESSIE electronic structure finite element code\cite{Kestyn:2020,NESSIE}, we have extracted two types of linear systems: (i) a spd matrix obtained after discretization of the Poisson equation, and (ii) complex symmetric matrices that arise from solving the Schr\"odinger eigenvalue problem using the FEAST solver \cite{PolizziFEAST,Gavin:2013,Polizzi2014,feast}.  
These matrices share the same dimensions, non-zero count, and sparsity structure. 
They are $796831 \times 796831$ matrices with $22751269$ elements. 
After the weighted spectral ordering is applied, these matrices have upper and lower half-band-widths of $35160$.  (the relative magnitude of the elements of these matrices is similar enough that the permutation generated by applying the weighted spectral ordering to one produces a reasonably good reordering for the other, so the weighted spectral ordering based on a matrix from (ii) has been used to reorder all of these matrices from the outset).
In addition, with only very tiny elements are present past a band width of $20000$, and so we limit the size of our off-diagonal submatrices $\M{B_i}$ and $\M{C_i}$ to $20k \times 20k$ elements. 

In the following sections, these subproblems will be used to evaluate the performance of \algname{}. 
For both of the problems we have here, the \algname{}-T configuration happens to be preferable, so this preconditioner will be the focus.

All of the measurements have been performed on a distributed system, consisting of compute nodes with the characteristics: 
\begin{itemize}
\item{CPU: 2$\times$ Intel(R) Xeon(R) Platinum 8352Y}
\item{Networking: Infiniband 2X HDR (100 Gb/sec)}
\end{itemize}
Each Xeon(R) 8352Y provides 32 cores, with two CPU packages per compute node, resulting in 64 cores per node.
Clock frequency may range from a 2.2GHz base freuency, to a 3.4GHz boost. 
In addition, AVX-512 vector extensions are availible. 
Intel's `oneAPI' suite of compilers and tools have been used to take full advantage of the hardware, 
\begin{itemize}
\item Intel Fortran Classic, 2022.0 release 
\item oneMKL, 2022.0 release (contains implementations of BLAS, LAPACK, as well as the Intel-MKL variant of PARDISO) 
\item Intel MPI, 2021.6 release
\end{itemize}

We will observe the performance improvement achieved by adding additional nodes, to increase the number of MPI processes. 
For all of the measurements, the systems have been configured such that four MPI processes are placed on each node, allocating $16$ cores to each MPI process (so, two MPI processes are bound to each CPU package). 
Multicore parallelism inside each MPI process is handled by Intel-MKL. 

For all configurations using preconditioners, a mixed precision strategy is used. 
The preconditiers applied are single-precision variants (single-real or single-complex as applicable). 
Most other computations are performed in double precision, except where noted otherwise. 
%Now, we may proceed to the measurements.\todo{I don't really like this transition, but I'll have to come back to it}

\subsubsection{NESSIE Poisson Matrix}\label{Poission}

\begin{figure}[htbp]
\caption{Comparison of \algname{}-T and Block Jacobi for preconditioning BiCGTtab iterations when solving NESSIE Poisson matrices. Markers indicate Conjugate Gradient or BiCGStab iterations.  }\label{Poisson_convergence_caption}
\centering 
  \subcaptionbox{Convergence behavior when using $12$ MPI processes \label{Poisson_convergence}}%
  [.49\linewidth]{\includegraphics[keepaspectratio,width=1.0\linewidth]{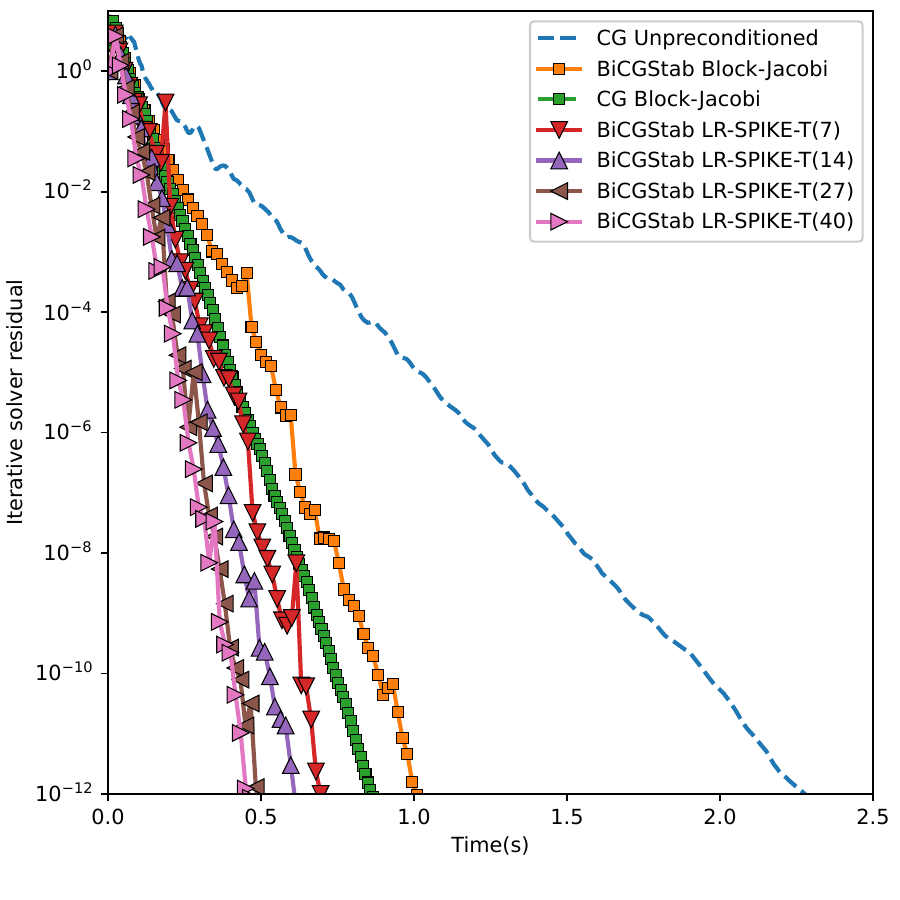}
  }%
  \hfill
  \subcaptionbox{Time, in seconds, required to reach convergence for each of the iterative solvers. \label{Poisson_stacked}}%
  [.49\linewidth]{\includegraphics[keepaspectratio,width=1.0\linewidth]{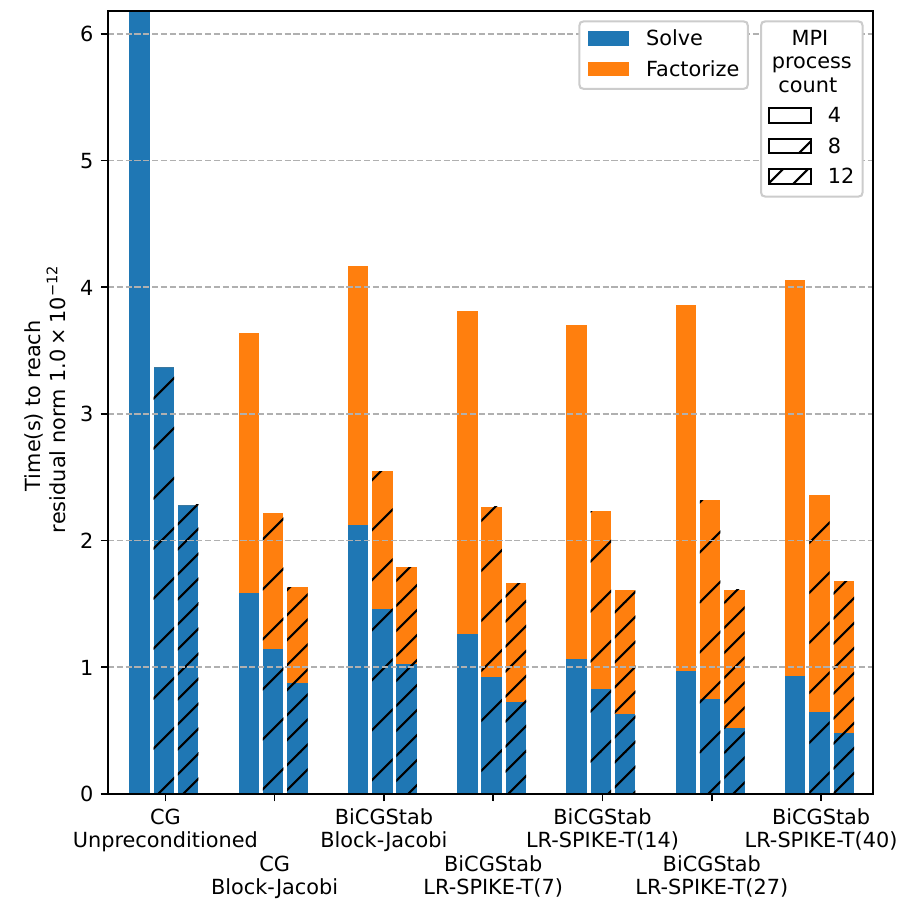}
  }%
  
\end{figure} 

\begin{table}[htbp]
\centering 
\begin{tabular}{|c|c|c|c|cccc|} \hline
Outer              & \tworow{Preconditioner} & \tworow{$n_{svd}$} & MPI           & Factorization         & Solve               & Combined            & Multiplication                \\       
      Solver       &                         &                    &     processes &               Time(s) &       Time(s)       &          Time(s)    &                 Count         \\\hline
                   &                         &                    &$4 $           &       N/A             &     $6.18$          &     $6.18 $         &       $697     $             \\    
                   & Unpreconditioned        &  N/A               &$8 $           &       N/A             &     $3.37$          &     $3.37 $         &       $697     $             \\    
\tworow{CG}        &                         &                    &$12$           &       N/A             &     $2.28$          &     $2.28 $         &       $697     $             \\\cline{2-8}
                   &                         &                    &$4 $           &$\bld{2.05}$           &     $1.59$          &$\bld{3.64}$         &       $80      $             \\    
                   & Block-Jacobi            &  N/A               &$8 $           &$\bld{1.07}$           &     $1.14$          &$\bld{2.21}$         &       $103     $             \\    
                   &                         &                    &$12$           &$\bld{0.76}$           &     $0.87$          &$\bld{1.63}$         &       $108     $             \\\hline   
                   &                         &                    &$4 $           &$\bld{2.04}$           &     $2.13$          &     $4.17 $         &       $53      $             \\    
                   & Block-Jacobi            &  N/A               &$8 $           &$\bld{1.09}$           &     $1.46$          &     $2.55 $         &       $64      $             \\    
                   &                         &                    &$12$           &$\bld{0.76}$           &     $1.03$          &     $1.79 $         &       $64      $             \\\cline{2-8}    
                   &                         &                    &$4 $           &     $2.55$            &     $1.27$          &     $3.81 $         &       $31      $             \\    
                   & \algname{}-T            &   7                &$8 $           &     $1.35$            &     $0.92$          &     $2.27 $         &       $41      $             \\    
                   &                         &                    &$12$           &     $0.93$            &     $0.73$          &     $1.66 $         &       $45      $             \\\cline{2-8}    
                   &                         &                    &$4 $           &     $2.63$            &     $1.07$          &     $3.70 $         &       $26      $             \\    
 BiCGStab          & \algname{}-T            &  14                &$8 $           &     $1.40$            &     $0.83$          &$\bld{2.23}$         &       $36      $             \\    
                   &                         &                    &$12$           &     $0.98$            &     $0.63$          &$\bld{1.61}$         &       $37      $             \\\cline{2-8}    
                   &                         &                    &$4 $           &     $2.89$            &     $0.97$          &     $3.86 $         &       $23      $             \\    
                   & \algname{}-T            &  27                &$8 $           &     $1.57$            &     $0.75$          &     $2.32 $         &       $32      $             \\    
                   &                         &                    &$12$           &     $1.09$            &     $0.52$          &$\bld{1.61}$         &       $30      $             \\\cline{2-8}    
                   &                         &                    &$4 $           &     $3.12$            &$\bld{0.93}$         &     $4.06 $         &  $\bld{22}     $             \\    
                   & \algname{}-T            &  40                &$8 $           &     $1.71$            &$\bld{0.65}$         &     $2.36 $         &  $\bld{27}     $             \\    
                   &                         &                    &$12$           &     $1.20$            &$\bld{0.48}$         &     $1.68 $         &  $\bld{28}     $      \\\hline     
\end{tabular}
\caption{Performance of \algname-T preconditioner, compared to Block Jacobi, for NESSIE Poisson type matrix. 
The best results, as well as those within $.02s$ of the best, have been marked with bold for each hardware configuration.  
Sparse multiplication counts are reported to allow easier comparisons between Conjugate Gradient (CG) and BiCGStab.
Convergence target is $\etothe{-12}$}.
\label{PoissonPerformanceTable}
\end{table}

In the NESSIE application, solve operations with this matrix must be done repeatedly and often involve only one right hand side (using a density functional theory and local exchange-correlations for first-principle electronic calculations \cite{Kestyn:2020}).
Since the factorization is performed only once, the cost of the preconditioner construction is quickly amortized. 
Additionally, high accuracy is required, and the matrix is  symmetric positive definite (SPD), and real. 

Because %the number of right hand sides is small, and 
the matrix is SPD, unpreconditioned Conjugate Gradient (CG) may be used to solve this matrix. 
Block Jacobi preconditioned CG is also a reasonable choice for distributed systems.  
For \algname{}-T, BiCGStab is used because the preconditioner is non-symmetric. 
Block Jacobi preconditioned BiCGStab is also considered, although it provides worse results %is strictly worse 
than CG for this matrix, it can be used as %to provide 
a more direct comparison to \algname{}-T. 
%As mentioned previously, \algname{}-T is used in all cases here. 

%Additionally, small $n_{svd}$ values were selected for \algname{}-T, because there is only one right-hand-side, and the matrix is not very challenging. 
%A larger value for $n_{svd}$ would result in a longer factorization time, which would not be sufficiently amortized across the already short solve times.
The specific values chosen for $n_{svd}$ were selected for best performance when computing the randomized SVD used to approximate the spikes. 
During this stage of the computation, solve operations must be performed upon right hand sides with dimensions $n_i \times \alpha n_{svd}$, where $\alpha$ is an oversampling factor used to improve the representation of the spikes in the randomized SVD.  
In this case, we use an oversampling factor of $\alpha = 1.2$.  
MKL PARDISO can exploit multicore parallelism across the solution vectors, and load balancing will be best if the number of right hand sides is an integer multiple of the number of threads used (particularly noticeable if the number of right hand sides is small).
Finally, the solve operations required for $\Vmat$ and $\Wmat$ may be handled simultaneously. 
As a result, values of $n_{svd}$ were selected to satisfy $2\lfloor 1.2 \times n_{svd} \rfloor=16$ ($16$ being the number of cores associated to each MPI process).

In \figs~\ref{Poisson_convergence} and ~\ref{Poisson_stacked}, parallel performance is shown. 
In terms of overall time, despite this matrix being SPD and fairly easy to solve, both 
preconditioners outperform unpreconditioned Conjugate Gradient substantially. 
Among the measurements using BiCGStab as the iterative solver, \algname{}-T outperforms Block Jacobi for overall computation time in every case (albeit sometimes by slim margins).

Block Jacobi preconditioned Conjugate Gradient and \algname{}-T preconditioned BiCGStab have similar overall computation times for most of the configurations.
%In particular, the overall times for \algname{}-T with $n_{svd}=14$ are almost identical to the more specialized solver.
But, unlike Block Jacobi, \algname{}-T provides the capability to increase the preconditioner quality at the cost of additional factorization time. 
This provides \algname{}-T a significant advantage in terms of solve time, with $n_{svd}=40$ ultimately achieving nearly a $2\times$ improvement over the alternative. 
%if desired, resulting in  \algname{}-T for this system. 
%prioritizing solve times over factorization. 
Performance and scalability improvements are expected to be even more significant for applications with multiple right-hand sides, such as those involving an exact-exchange term (including Hartree-Fock) in electronic structure calculations \cite{Braun:2014}.

%For this matrix, all of the solvers scale well as additional partitions are added, as can be seen in \fig~\ref{Poisson_stacked}.
As can be seen in \fig~\ref{Poisson_stacked}, all of the solvers improve in performance as additional MPI processes are added. 
This is not a given, as seen in section~\ref{Sherman5_subsec} the preconditioner may lose some effectiveness as the number of partitions is increased. 
In this case, as can be seen in table~\ref{PoissonPerformanceTable}, there is a slight increase in the number of iterations (measured in terms of sparse multiplication operations) as the number of partitions is increased. 
However, the rise in iteration count remains minimal, ensuring good scalability.
Next, we will look at a more challenging problem, for which the quality of the \algname{} preconditioner becomes essential.

\subsubsection{FEAST using \algname{}}

%The first set of measurements come from the use of FEAST to solve a general eigenvalue problem in the NESSIE application. 

Recent development of the FEAST eigenvalue solver has created an opportunity to use a fast approximate linear system solver \cite{Gavin2018,Li:2025}.
This role can be filled by an iterative solver given a loose convergence criteria, or by a direct solver working in reduced precision. 
As we will see shortly, \algname{} is quite suitable for this application.

For the specific eigenvalue problem discussed here, the goal is to find a subset of the negative eigenvalues near zero for a generalized, Hermitian eigenvalue problem. 
\begin{align}
  \label{EIG-beginning}
  \M{A} \M{X} =\M{B} \M{X} \M{\Lambda} , 
\end{align}
where $\M{A}$ and $\M{B}$ are sparse Hermitian matrices. 
In addition, $\M{B}$ is positive definite, though $\M{A}$ is not. 
Finally, we will assume there are $m$ eigenvalues of interest, and they lay in some range $\left(\lambda_{min}, \lambda_{max}\right)$

FEAST computes eigenvalues via filtered subspace iteration.
The filter used is an approximation of the spectral projector, created by approximating a contour integration along a path in the complex plane around the eigenvalues of interest. 
The operator to implement is
 \begin{align}
	\rho(\M{\inv{B}}\M{A}) = \frac{1}{2\pi i} \oint\limits_{C} \inv{(z\M{B} - \M{A})}\M{B} dz \approx \sum_{j=1}^{n_{c}} \omega_{j} \inv{(z_j \M{B} - \M{A})}, \label{contour_intergration}
 \end{align}
where $z_j$ and $\omega_{j}$ are the points and corresponding weights used to perform the numerical integration around a path, $C$, enclosing the eigenvalues of interest. 
Eigenvalues of \eqn~\ref{EIG-beginning} outside the contour correspond to eigenvalues of \eqn~\ref{contour_intergration} with low magnitude, and so, given a large enough subspace, the subspace iteration will converge on those contained within the contour (if the integration could be performed exactly, rather than via a numerical approximation, the eigenvalues of \eqn~\ref{EIG-beginning} outside of the contour would correspond to zero eigenvalues of \eqn~\ref{contour_intergration}, and FEAST would converge in one step).
The primary computational cost of FEAST comes from solving the shifted linear systems corresponding to the contour points, $\inv{(z_j \M{B} - \M{A})}$.

Increasing the number of contour points, $n_c$, will improve the quality of this numerical integration. 
This produces a better filter for the subspace iteration, resulting in faster convergence for the eigenvalue problem.
Because the linear systems associated with the contour points may be solved simultaneously, this provides a source of parallelism. 

Originally FEAST used direct solvers (such as PARDISO) to perform the linear system solve operations for the contour points.
However, recent FEAST developments have improved support for iterative solvers.
In particular, a residual inverse iteration based formulation has been implemented \cite{Gavin2:2018,feast}, which allows for convergence of the eigenvalue problem even if the linear system is not solved exactly (originally this was intended for solving nonlinear eigenvalue problems, but it is applicable to more conventional linear eigenvalue problems as well).
With this new formulation, the contour point linear system solve operations may be performed to a loose tolerance.
This will introduce a new limit on the rate of convergence for the subspace iterations, equal to the tolerance to which the linear systems are solved. 
%However, as long as the linear system solve operations are performed to a tolerance tighter than the convergence rate imposed by the numerical approximation of the contour integral, this has not been observed to impact the overall convergence of the eigenvalue problem. 
However, when the linear system solve operations are performed accurately enough to meet the convergence rate imposed by the numerical approximation of the contour integral, it can be shown that the overall convergence rate of the eigenvalue problem is not harmed \cite{feast}.
%However, when the linear system solve operations are performed to a tolerance equal to the convergence rate imposed by the numerical approximation of the contour integral, this has not been observed to impact the overall convergence of the eigenvalue problem. 
%\todo{Maybe a note about the version of FEAST in which this will be available to the general public, or where more detail can be found}

%A full analysis of the impact of inexact linear system solve operations inside FEAST on the convergence rate of the eigenvalue problem has not been perfomed. 

For this particular application, we seek a subset of the eigenvalues which lie in the range $\left(\lambda_{min}=-65,\lambda_{max}=-4\right)$. 
This range of eigenvalues covers the energy states of the valence electrons for the 40zgnr nanostructure inspected by this NESSIE simulation \cite{Polizzi:2015,Kestyn:2020}.
These are not extremal eigenvalues, there is an additional cluster of eigenvalues around -275 corresponding to the core electrons of the system, and there are a great number of eigenvalues above -4.0, extending past zero into the right side of the complex plane, corresponding to the continuum spectrum. 
FEAST also requires an over-estimate of the number of eigenvalues in the range of interest, to size the search subspace. 
It is desirable for this value, $m_0$, to be in the range of $1.5m$ to $2m$, where $m$ is the number of eigenvalues in the range of interest. 
An estimate of $m_0=700$ is provided by an earlier stage of the simulation (and in fact there are $m=379$ eigenvalues in this range).
A visualization of the problem can be seen in \fig{}~\ref{40zgnr_contour}.

 \begin{figure}[h]
 \centering 
 \caption{The FEAST eigenvalue problem to be solved, plotted on the complex plane. The eigenvalues of interest, labeled `valence eigenvalues' correspond to the energy levels of valence electrons. Because the problem is symmetric, the eigenvalues have no imaginary part. Each FEAST contour point produces a linear system problem with a complex shift, $\M{M}=z\M{B}-\M{A}$. Contour points which are near eigenvalues correspond to more difficult linear system problems, as the system to be solved becomes closer to singular. However, and the existence of a wide band-gap provides some flexibility in defining the contour, which allows for the selection of contour points with good separation from the eigenvalues.}\label{40zgnr_contour}
 \includegraphics[keepaspectratio,width=.75\linewidth]{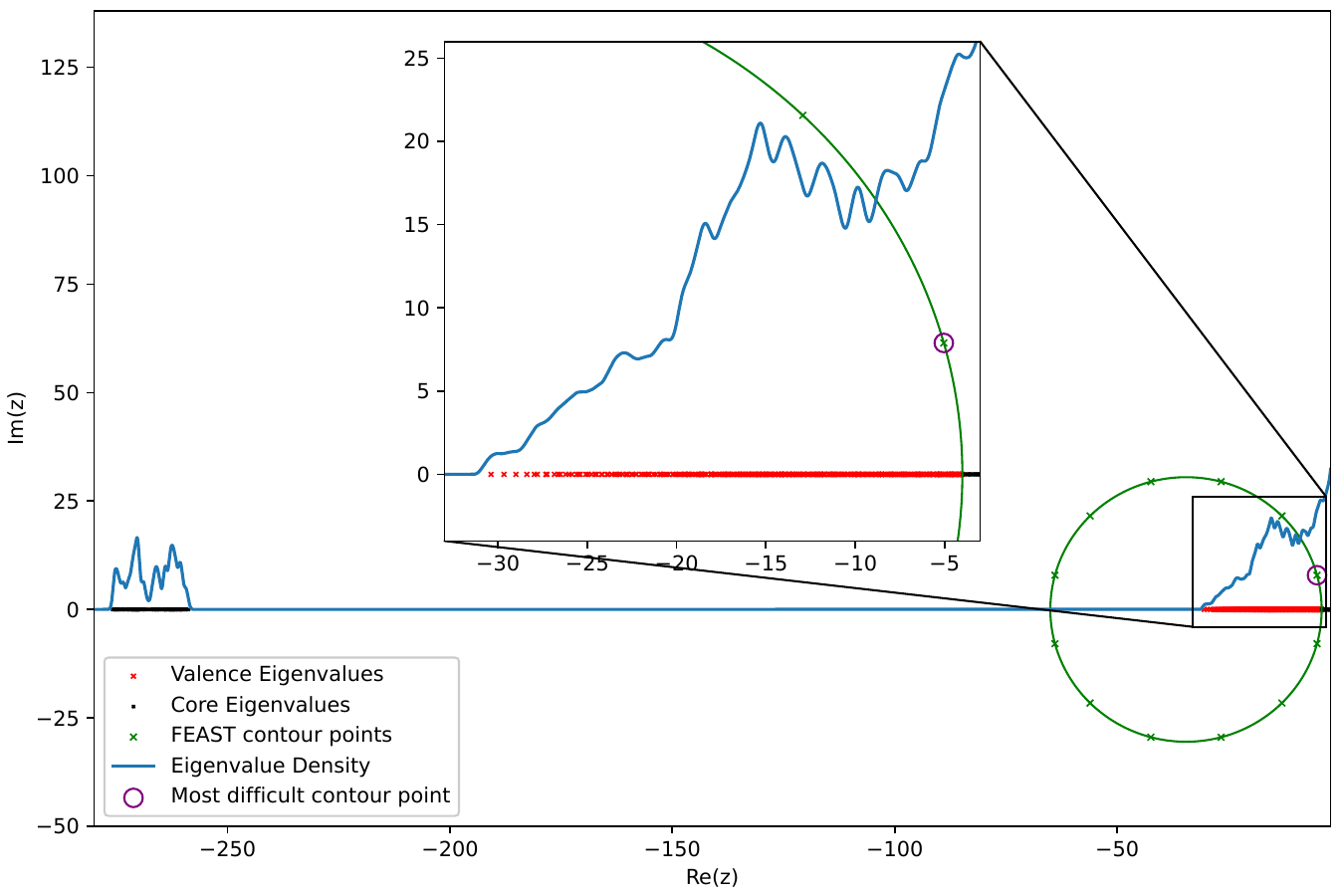}
 \end{figure} 

A total of twelve contour points are shown in \fig{}~\ref{40zgnr_contour}, which provides six linear systems to solve (for symmetrical problems, FEAST internally exploits symmetry to handle the contour points with negative imaginary parts).
The linear systems corresponding to contour points near a large numbers of eigenvalues tend to be the most difficult to solve.
%(for intuition sake, note that in the unlikely event that a contour point were placed exactly on an eigenvalue, the corresponding linear system would be singular). 
In this case the most troubling contour point is the one farthest to the right, marked with a purple circle in \fig{}~\ref{40zgnr_contour}. 
Consequently, we will begin by looking at this linear system in isolation to gain some insight into the performance of the solver. % and then use this information to select  values for $n_{svd}$.  
For these runs the number of right hand sides has been reduced to $n_{rhs}=32$. 
Ultimately inside FEAST, the number of right hand sides used when solving will be $n_{rhs} = m_0 = 700$, but $32$ is sufficient to characterize the performance of the solver.

%%% some changes here
%The behavior of the solvers is shown in ~\fig{}\ref{40zgnr_performance}.
%Two convergence thresholds are considered. 
The convergence behavior and scalability performance of the solvers are shown in ~\fig{}\ref{40zgnr_performance}.
\fig{}\ref{Hardshift_convergence} clearly shows that \algname{}-T outperforms Block Jacobi for this matrix using $12$ MPI processes/partitions. 
For these runs, two convergence thresholds are considered.
A threshold of $\etothe{-12}$ is used to represent convergence to around working precision. 
In addition, time to reach  $1.5 \times \etothe{-1}$ is shown. 
The eigenvalue problem from which this contour point has been extracted will require 15 iterations to converge to $\etothe{-12}$ when using the previously described FEAST configuration, even if the linear systems are solved to machine precision.
This is the convergence rate imposed by the quality of the contour integration, and so a linear system solve convergence threshold tighter than $\etothe{-12/15} \approx 1.5 \times \etothe{-1}$ is wasteful. A method of discovering of this limit is not discussed here.
%% this added here from below
Of particular interest are $n_{svd}=512$ and $n_{svd}=768$, which converge to $1.5 \times \etothe{-1}$ in just two and one BiCGStab iterations respectively.

%\textemdash because multiple eigenvalue problems with similar matrices must be solved for the NESSIE application, it is assumed that this could be discovered at an earlier stage in the computation (alternatively, dynamically adjusting the linear system solve tolerance based on the rate of convergence of the eigenvalue problem may be a possible extension to FEAST).  
 \begin{figure}[h]
 \centering 
 \caption{
Performance characteristics of BiCGStab, when preconditioned by Block Jacobi or \algname{}-T, for the linear system solve operation corresponding to the most difficult FEAST contour point shown in ~\fig \ref{40zgnr_contour}. 
$32$ right hand sides are used. 
\textbf{(a)} Convergence behavior when using $12$ MPI processes. 
Markers indicate BiCGStab iterations. 
Inset, the convergence behavior down to the `loose' convergence which FEAST can tolerate for this problem, $1.5 \times \etothe{-1}$, is shown. 
\textbf{(b)} Parallel performance for each of the solvers. 
\textbf{(b)}-top shows the time, in seconds, to reach a `tight' convergence of $\etothe{-12}$ including the factorization time, while \textbf{(b)}-bottom shows just the solve time to reach the loose convergence needed by FEAST, $1.5 \times \etothe{-1}$.
}\label{40zgnr_performance}
  %\subcaptionbox{BiCGStab residual time. Markers indicate individual iterations.}%
  \subcaptionbox{\label{Hardshift_convergence}}%
  [.5\linewidth]{
  \includegraphics[keepaspectratio,width=\linewidth]{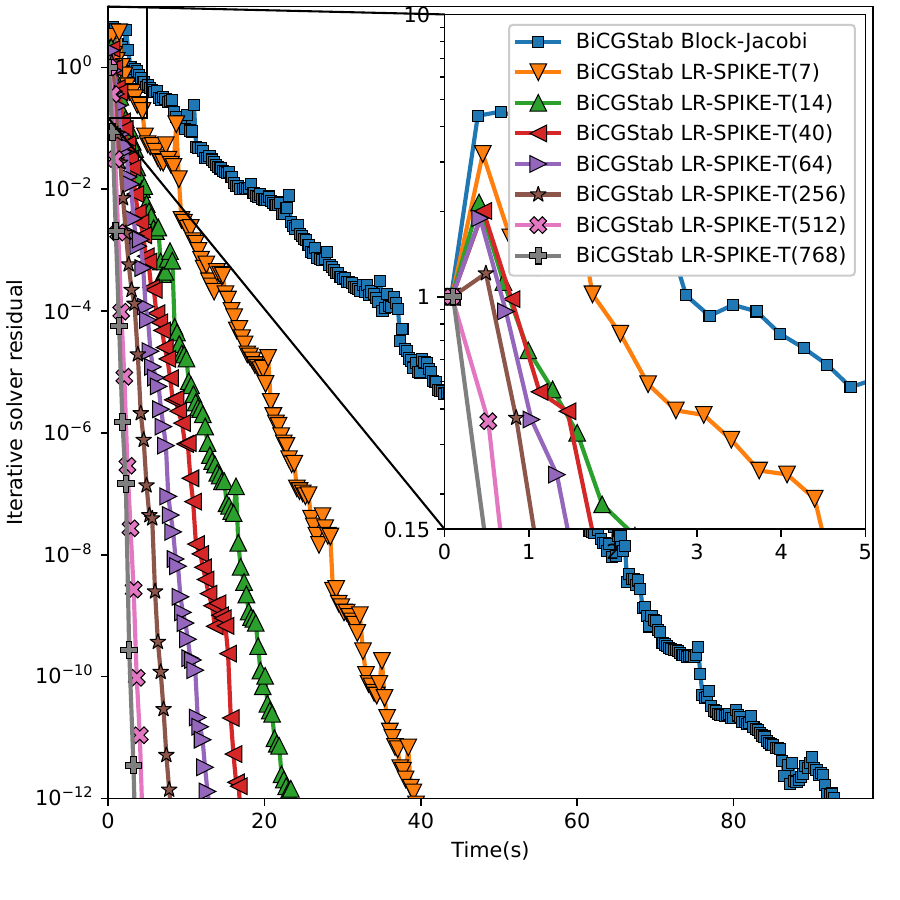}
  }%
  \subcaptionbox{\label{Both_stacks}}%
  [.5\linewidth]{
  \begin{tabular}{c}
  \includegraphics[keepaspectratio,width=\linewidth]{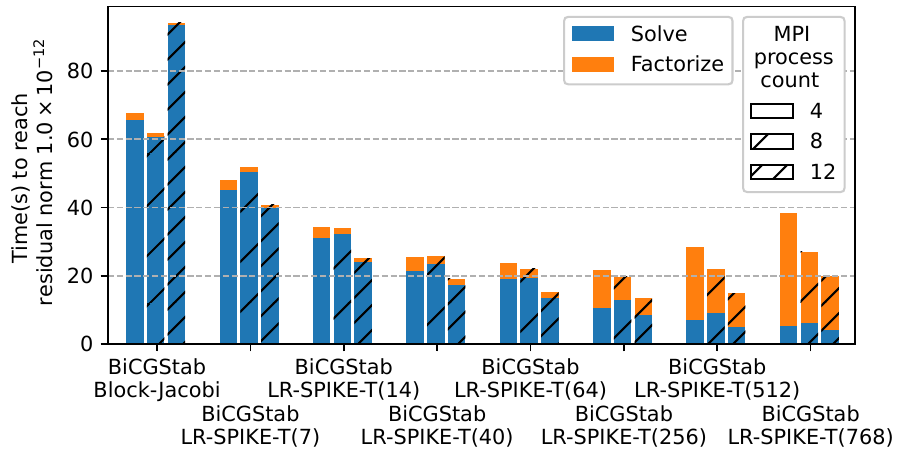} \\ 
  \includegraphics[keepaspectratio,width=\linewidth]{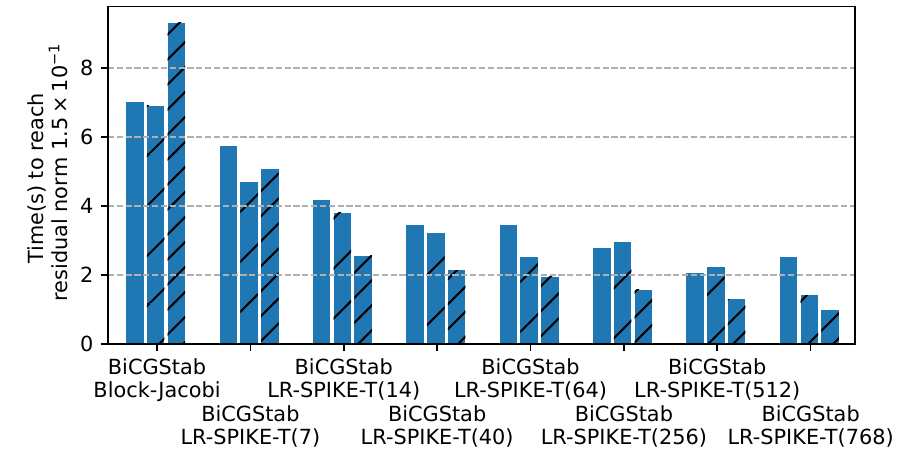}
  \end{tabular}%
  }
 \end{figure} 
%%% changes (added before)
%For this matrix, \algname{}-T outperforms Block Jacobi substantially. 

%% changes
%In \fig{}~\ref{Both_stacks}, the time to reach the two different convergence thresholds is shown.
\fig{}~\ref{Both_stacks} show the performance scalability  from 4 to 12 MPI/partitions for the total, factorization and solve stages, as well as the time to reach the two different convergence thresholds.
It can be see in \fig{}~\ref{Both_stacks}-top that even with just 32 right hand sides, the additional cost of the factorization is easily amortized for most cases. 
Additionally, in cases where the factorization time is still significant, $n_{svd}=256$ and $n_{svd}=512$, the factorization time is quite scalable. 
A noticeable loss of performance for the Block Jacobi preconditioner can be seen the number of partitions is increased.
\algname{} mitigates this limitation for small values of $n_{svd}$, less than 40,  and essentially removes it for the larger values.  

%% I added some of it before
%In \fig{}~\ref{Hardshift_convergence}, convergence behavior of the BiCGStab can be seen, for $12$ MPI processes. 
%Inset, the behavior for convergence to the looser FEAST-appropriate tolerance is shown.
%Of particular interest are $n_{svd}=512$ and $n_{svd}=768$, which converge to $1.5 \times \etothe{-1}$ in just two and one BiCGStab iterations respectively. 

\algname{} with $n_{svd}=7$ or $14$ is an improvement to Block Jacobi, and the additional factorization time is quite minimal. 
However, due to the presence of many more right hand sides when performing the solve operations inside FEAST, as well as the multiple outside iterations, solve time is a much higher priority in this application. 
%% I added the first sentence
\fig{}~\ref{Both_stacks}-bottom shows the scalability of the solve time to reach the loose convergence needed by FEAST. 
In the following, \algname{} preconditioners with values of $n_{svd} \geq 40$ were selected from among those tested here for use in the following run of FEAST.

In table~\ref{FeastPerformanceTable}, the result of using these solvers inside FEAST is shown.
Times are reported specifically for the linear system solvers, as well as the overall time consumed by the computation.
As mentioned previously, the linear systems associated with the contour points vary in difficulty, with contour points near the real line and deeper in the eigenvalue spectrum producing more difficult linear systems. 
Because FEAST handles the these linear system simultaneously, the overall computation will have to wait for whichever linear system takes the longest to solve. 
In the case of the iterative solvers this should be the most difficult point, which will require additional BiCGStab iterations to converge.
As such, the data collected for the solve operations is focused on ``most difficult'' linear system of the set (as found simply by measuring all of them and selecting the worst results for each entry). 

In this application \algname{}-T outperforms Block Jacobi substantially, in some cases by an order of magnitude. 
This is unsurprising, given the nature of the problem\textemdash the overall computation time is mostly dominated by the cost of solving linear systems for this system, giving \algname{} a natural advantage.
Both preconditioners do lose some effectiveness in the transition from four to eight partitions. 
But \algname{}-T does not further degrade when the number of partitions is increased to twelve, producing a net gain in overall performance in most cases (also note that, because four MPI processes have been placed on each compute node, the transition from eight to twelve MPI processes is more indicative of the scalability of the solver than the transition from four to eight).  %%%% Question here

 \begin{table}[htbp]
 \centering
 \caption{Performance measurements for FEAST eigenvalue solver to find $m=379$ eigenpairs (subspace size $m_0 = 700$) for a system size $n \approx 800k$, using various linear system solvers for the contour point solve operations. 
 Six contour points are used in all cases, resulting in a total number of MPI processes which is equal to the total number of partitions times six. 
 Each MPI process uses 16 cores for Intel-MKL parallel threading (a total of $12\times 6\times 16=1152$ cores are then used for the largest runs).
 The factorization of the linear system at each contour points is performed only once, in the first FEAST iteration.  
 The best timings in each field for each hardware configuration are marked in bold. 
 Convergence target for the eigenvalue problem is $\etothe{-12}$.
 For BiCGStab linear system solve convergence, a target of $1.5 \times \etothe{-1}$ has been selected as described in the text.
 $\dagger{}$ indicates measurements which vary from contour point to contour point. 
 For those measurements, the `worst' (slowest or most iterations) result among all contour points is shown, as that becomes the bottleneck.
 ``Prec calls'' reports a count of the total number of calls to the preconditioner or direct solver, for the most difficult contour point, aggregated over all FEAST iterations. The total FEAST time includes not only the total time spent solving the linear systems, but also the FEAST overhead: communication and other auxiliary computations.}
\label{FeastPerformanceTable}                                                                                                                                                                                    
\begin{tabular}{|c|c|c|cc|cc|cc|c|} \hline                                                                                                                                                                           
    \tworow{ Iterative }    &                          &                 &                    & \tworow{Total}  &\tworow{\#FEAST} & \tworow{\#Prec.            }      &  \multicolumn{2}{c|}{Aggregate                }  & Total      \\ 
    \tworow{  Solver   }    &   Preconditioner         &  $ n_{svd} $    &     \#Partition    & \tworow{\#MPI}  &\tworow{Iters} & \tworow{Calls$ ^\dagger $}          &  \multicolumn{2}{c|}{          solver times(s)}  & FEAST      \\  
                            &                          &                 &                    &                 &      &                                     &    Fac. $ ^\dagger $  & Sol. $ ^\dagger $       & Time(s)    \\\cline{1-10}                    
                                                                                                                                                                                                                                                       
                            &                          &                 &  $4 $&$ 24 $  &  $15  $  &  $526   $    &  $\bld{2.2}   $  &  $ 1969.9 $   &  $2089.3$          \\                           
                            &Block Jacobi              &    N/A          &  $8 $&$ 48 $  &  $16  $  &  $955   $    &  $\bld{1.2}   $  &  $ 2500.7 $   &  $2627.9$          \\                           
                            &                          &                 &  $12$&$ 72 $  &  $15  $  &  $2011  $    &  $\bld{0.9}   $  &  $ 3746.2 $   &  $3921.0$          \\\cline{2-10}
                                                                                                                                                               
                            &                          &                 &  $4 $&$ 24  $  &  $15  $  &  $97    $    &  $    3.9     $  &  $  444.4 $  &  $581.5 $          \\                           
                            &                          &    40           &  $8 $&$ 48  $  &  $16  $  &  $165   $    &  $    2.1     $  &  $  535.6 $  &  $599.4 $          \\                           
                            &                          &                 &  $12$&$ 72  $  &  $15  $  &  $194   $    &  $    1.6     $  &  $  444.7 $  &  $487.4 $          \\\cline{3-10}
                                                                                                                                                         
   \tworow{BiCGStab}        &                          &                 &  $4 $&$ 24  $  &  $15  $  &  $78    $    &  $    4.5     $  &  $  372.8 $  &  $510.7 $          \\                           
\tworow{$ 1.5 \times \etothe{-1} $ }&                  &    64           &  $8 $&$ 48  $  &  $16  $  &  $148   $    &  $    2.5     $  &  $  484.0 $  &  $549.7 $          \\                           
                            &  \tworow{\algname{}-T}   &                 &  $12$&$ 72  $  &  $15  $  &  $163   $    &  $    1.8     $  &  $  389.3 $  &  $432.1 $          \\\cline{3-10}  
                                                                                                                                                         
                            &                          &                 &  $4 $&$ 24  $  &  $15  $  &  $31    $    &  $   10.3     $  &  $ 227.5  $  &  $341.8 $          \\                           
                            &                          &   256           &  $8 $&$ 48  $  &  $15  $  &  $93    $    &  $    6.3     $  &  $ 342.6  $  &  $405.0 $          \\                           
                            &                          &                 &  $12$&$ 72  $  &  $15  $  &  $90    $    &  $    4.7     $  &  $ 245.8  $  &  $289.2 $          \\\cline{3-10}                
                                                                                                                                                         
                            &                          &                 &  $4 $&$ 24  $  &  $15  $  &  $16    $    &  $   19.4     $  &  $ 159.3  $  &  $276.2 $          \\                           
                            &                          &   512           &  $8 $&$ 48  $  &  $15  $  &  $37    $    &  $   12.3     $  &  $ 172.5  $  &  $246.0 $          \\                           
                            &                          &                 &  $12$&$ 72  $  &  $15  $  &  $31    $    &  $    9.4     $  &  $ 109.9  $  &  $158.6 $          \\\cline{1-10}                
                                                                                                                                                         
                            &   MKL-PARDISO            &   N/A           &  $1 $&$ 6   $  &  $15  $  &  $16    $    &  $   15.0     $  &  $   132.8$  &  $  589.3    $     \\ \cline{2-10}                         
                                                                                  
                            &                          &                 &  $4 $&$ 24  $  &  $15  $  &  $16    $    &  $   19.5     $  & $\bld{75.8}$ &  $\bld{191.7}$     \\  
                            &                          &   512           &  $8 $&$ 48  $  &  $27  $  &  $28    $    &  $   12.3     $  &  $    72.2$  &  $183.9 $          \\ 
                            & \tworow{\algname{}-T}    &                 &  $12$&$ 72  $  &  $21  $  &  $22    $    &  $    9.7     $  &  $    43.0$  &  $103.3 $          \\\cline{3-10} 
                                                                                                                                                         
                            &                          &                 &  $4 $&$ 24  $  &  $15  $  &  $16    $    &  $   31.5     $  &  $   79.1$   &  $     207.4 $     \\                           
   None                     &                          &   768           &  $8 $&$ 48  $  &  $18  $  &  $19    $    &  $   19.6     $  &$\bld{53.1}$  &  $\bld{140.5}$     \\  
                            &                          &                 &  $12$&$ 72  $  &  $15  $  &  $16    $    &  $   15.1     $  &$\bld{34.8}$  &  $\bld{ 87.0}$     \\\cline{2-10}                
                                                                                                                                                         
                            & \tworow{MKL Cluster}     &                 &  $4 $&$ 24  $  &  $15  $  &  $16    $    &  $   13.1     $  &  $   159.7$  &  $  273.2    $     \\                           
                            & \tworow{Sparse Solver}   &  N/A            &  $8 $&$ 48  $  &  $15  $  &  $16    $    &  $   12.9     $  &  $  1831.5$  &  $ 1899.7    $     \\                           
                            &                          &                 &  $12$&$ 72  $  &  $15  $  &  $16    $    &  $   13.0     $  &  $   206.3$  &  $  255.4    $     \\\cline{1-10}                
                                                                                                                                                                                       
\end{tabular}   
\end{table}

It was also noted in \fig{}~\ref{40zgnr_performance} that, with sufficiently large $n_{svd}$, \algname{}-T might not require the assistance of BiCGStab to attain the accuracy necessary for FEAST. 
These configurations are shown in the lower rows of the Table~\ref{FeastPerformanceTable}. 
\algname{}-T could be thought of as a very approximate direct solver in this configuration, and so the default direct solvers used in FEAST, MKL Pardiso (shared memory) and the Intel-MKL Cluster Sparse Solver (distributed), are used as a points of comparison. 
We have made our best effort to configure these solvers appropriately for this problem: iterative refinement has been disabled (as FEAST does not require high accuracy), and single-precision complex data types have been used. 

The Cluster Sparse Solver performs well when used with four MPI processes, but when the computation is spread across additional compute nodes, performance degrades rather than improving.
We were unable to locate the cause of the performance degradation for eight MPI processes in the case of the Cluster Sparse Solver, but believe the results from four and twelve MPI processes are more representative of the general behavior of the solver (the performance degradation appears to be the result of an unlucky pattern in the fill-in causing heavy communication between the nodes). 

For Intel-MKL Pardiso, a single MPI process is associated with each FEAST contour point. 
For this case, it is notable that the total aggregate Intel-MKL Pardiso time ends up being much smaller than the total FEAST simulation time ($589.3$s). 
The latter includes communication costs due to MPI reduction operations across all of the contour points while performing the numerical integration. 
As the number of partitions grows, this overhead time diminishes linearly, a trend consistently observed across all the reported distributed solvers.

When good values are selected for $n_{svd}$, \algname{}-T provides smoother scaling and better overall performance than the Intel-MKL Cluster Sparse solver for this problem. 
With a value of $n_{svd}=768$, for four and twelve partitions, the \algname{}-T solver is sufficiently accurate to match the FEAST contour integration, resulting in no additional FEAST iterations.
For eight partitions, \algname{}-T does see a slight increase in the number of FEAST iterations, although not enough to seriously harm overall performance. 
It should be emphasized that, while $768$ is large relative to the other values discussed in this section for $n_{svd}$, it is still quite small relative to the actual problem, around $2\%$ of the matrix half band width. 

\algname{}-T with $n_{svd}=512$  does not meet the desired linear system accuracy for FEAST for more than four \spike{} partitions, when used without BiCGStab.
However, an interesting comparison can be made between the two cases where \algname{}-T is used with $n_{svd}=512$\textemdash as a preconditioner for BiCGStab, and as an approximate direct solver.  
Using the twelve partition case as an example, \algname{}-T used as a preconditioner for BiCGStab requires two applications per FEAST iteration (in other words, a single full BiCGStab iteration) to reach the linear system convergence target.
On the other hand, when \algname{}-T is used without BiCGStab, the ideal linear system solver accuracy for FEAST is not quite met, resulting in additional FEAST iterations.
This presents a trade-off, we may accept a somewhat increased number of FEAST iterations, in exchange for a significantly less computationally expensive linear system solve operation. 
In this case, because the FEAST iterations are not very costly and including BiCGStab doubles the number of calls to \algname{}-T, it is preferable to use \algname{}-T as an inaccurate direct solver. 

% \subsubsection{Poisson}
\section{Conclusion}

The \spike{} framework for crafting parallel linear system solvers has been extended to include a truncated SVD based low-rank approximation of the eponymous ``spikes.'' 
This new scheme allows the user to accept some inaccuracy, in exchange for a more compact representation of the spikes. 
An efficient randomized SVD process is shown to be effective for generating these approximate spikes.  
Subsequently, they are used to construct approximate versions of classic \spike{} building-blocks, which may be dropped into the \spike{} framework to create new approximate variants of the solvers, suitable for use as preconditioners and, in some cases, even usable as very approximate direct solvers. 

Three of these low-rank-approximation enhanced solvers have been implemented and discussed here: 
\algname{}-I, which approximates the full spikes using their low-rank equivalents; \algname{}-T a low-rank enhanced version of the truncated \spike{} solver, particularly effective for diagonally dominant matrices; and \algname{}-OTF, a low-rank-enhanced variant of the classic ``\spike{}-on the fly'' solver. 
The new solvers are tested on a selection of publicly available matrices from SuiteSparse, and are compared in terms of effectiveness against some popular ``black-box'' preconditioners (which is to say, preconditioners which work fairly well for a wide range of problems, and do not require extensive study of the underlying matrix structure).
The \algname{} solvers are shown to be generally much more effective than the Block Jacobi preconditioner, and competitive with incomplete LU, for reducing the number of iterations required by BiCGStab. 
In addition, the performance and scalability of \algname{}-T have been shown using a problem extracted from the NESSIE, our finite element electronic structure code. 
A symmetric positive definite matrix problem has been extracted from the simulation. 
In this scenario, \algname{}-T performs well using very low rank approximations of the spikes.
In terms of overall performance, it compares favorably to Block Jacobi preconditioned BiCGStab. % and to unpreconditioned Conjugate Gradient.
It is also competitive with a solver which takes advantage of the SPD nature of this problem, Block Jacobi preconditioned Conjugate Gradient (roughly equal in terms of overall combined factorization and solve time, with \algname{}-T performing significantly better in the solve stage).

Finally, \algname{}-T is shown in use with the FEAST eigenvalue solver, which requires a performant solver for shifted linear system problems. 
From the NESSIE simulation a generalized, symmetric eigenvalue problem is extracted. 
In this example, non-extremal eigenvalues are sought (which is to say, eigenvalues not close to the edge of the spectrum, resulting in relatively difficult shifted linear system problems). 
\algname{}-T, when used as a preconditioner, drastically outperforms Block Jacobi for this problem. 
Furthermore, a particular synergy between FEAST and \algname{} is noted here; FEAST does not demand that the shifted linear system problem be solved exactly to working precision, and because the \algname{} low-rank approximation can become fairly close to an accurate representation of the spikes. 
As a result it is possible to use \algname{}-T as a ``low-accuracy direct solver,'' instead of as a preconditioner for BiCGStab. 
In this configuration, \algname{}-T is shown to outperform high quality direct solvers\textemdash the MKL variant of PARDISO, and the MKL Cluster Sparse Solver.  

\bmsection*{Acknowledgments}
This work utilized resources from Unity, a collaborative, multi-institutional high-performance computing cluster managed by UMass Amherst Research Computing and Data.

\bibliography{spike_svd}

%\bmsection*{Author Biography}

%\begin{biography}{\includegraphics[width=76pt,height=76pt,draft]{empty}}{
%{\textbf{Author Name.} Please check with the journal's author guidelines whether
%author biographies are required. They are usually only included for
%review-type articles, and typically require photos and brief
%biographies for each author.}}
%\end{biography}

\end{document}